\documentclass{article}[12pt]

\newtheorem{fed}{\textbf{Definition}}[section]
\newtheorem{thm}[fed]{\textbf{Theorem}}
\newtheorem{lemma}[fed]{\textbf{Lemma}}

\newtheorem{rem}[fed]{\textbf{Remark}}
\newtheorem{prop}[fed]{\textbf{Proposition}}
\newtheorem{cor}[fed]{\textbf{Corollary}}
\usepackage{amssymb,graphicx,epsfig,psfrag,epic,eepic,latexsym}
\usepackage{amsmath}
\usepackage{mathrsfs}
\usepackage{xcolor}
\usepackage[all]{xy}
\def\Nabla#1{\nabla\kern-.5ex{}_{#1}}

\newcommand{\C}{\mathcal{C}}

\newcommand{\N}{\mathbb{N}}

\newcommand{\Z}{\mathbb{Z}}
\newcommand{\R}{\mathbb{R}}

\newcommand{\ddt}{\frac{d}{dt}}
\newcommand{\dds}{\frac{d}{ds}}

\usepackage{verbatim}

\newcommand{\BB}{\mathcal{B}}

\newcommand{\HH}{\mathcal{H}}

\newcommand{\LL}{\mathcal{L}}

\newcommand{\QQ}{\mathcal{Q}}

\newcommand{\FF}{\mathcal{F}}
\newcommand{\la}{\langle}
\newcommand{\ra}{\rangle}
\newcommand{\wt}{\widetilde}

\newcommand{\eps}{\varepsilon}
\newcommand{\into}{\hookrightarrow}
\newcommand{\ol}{\overline}

\newcommand{\im}{{\rm im}\,}
\newcommand{\coker}{{\rm coker}\,}
\newcommand{\Crit}{{\rm Crit}}
\newcommand{\Func}{{\rm Func}}
\newcommand{\ind}{{\rm ind}}

\begin{document}
\title{Nondegeneracy and integral count of frozen planet orbits in helium}
\author{Kai Cieliebak, Urs Frauenfelder, Evgeny Volkov}
\maketitle
\begin{abstract}
%In~\cite{frauenfelder1} it was shown that for the helium atom with
%mean interaction between the electrons, there exists a unique frozen
%planet orbit. In this paper we improve this result by showing
%that the unique frozen planet orbit is nondegenerate and . The main
%idea of this paper is to
We study a family of action functionals whose critical points
interpolate between frozen planet orbits for the helium atom with
mean interaction between the electrons and the free fall. The rather
surprising first result of this paper asserts that for the whole
family, critical points are always nondegenerate. This implies that
the frozen planet orbit with mean interaction is nondegenerate and
gives a new proof of its uniqueness.
%Moreover, we explain that the moduli space of critical points
%of these action functionals is compact. Since there is a unique free fall orbit, this
%implies uniqueness of the frozen planet orbit. 
As an application, we show that the integral count of frozen planet
orbits with instantaneous interaction equals one. For this, we prove
orientability of the determinant line bundle over the space of
self-adjoint Fredholm operators with spectrum bounded from below, and
use it to define an integer valued Euler characteristic for Fredholm
sections whose linearization belongs to this class.  
\end{abstract}

%\tableofcontents

%%%%%%%%%%%%%%%%%%%%%%%%%%%%%%%%%%%%%%%%%%%%%%%%%%%%%%%%%%%%%%%%
\section{Introduction}\label{sec:intro}
%%%%%%%%%%%%%%%%%%%%%%%%%%%%%%%%%%%%%%%%%%%%%%%%%%%%%%%%%%%%%%%%

{\em Frozen planet orbits} are periodic orbits in the helium atom
in which both electrons move on a line on the same side of the
nucleus. The inner electron undergoes consecutive collisions with the
nucleus, while the outer electron (the actual ``frozen planet'') remains
almost stationary at some distance.
See~\cite{tanner-richter-rost, wintgen-richter-tanner} for numerical
%marginpar{\tiny Check references.}
evidence for such orbits and a discussion of their role in the
semiclassical treatment of the helium atom. 

When trying to prove the existence of frozen planet orbits, one faces the
difficulty that they cannot be obtained as perturbations of the system
without interaction between the electrons. In order to
deal with this problem, the second author replaced in~\cite{frauenfelder1} 
the instantaneous interaction between the two electrons by a mean
interaction and proved that in this case for every negative energy
there exists a unique frozen planet orbit. 
Building on work of Barutello, Ortega and Verzini~\cite{barutello-ortega-verzini},
we introduced in~\cite{cieliebak-frauenfelder-volkov}
two functionals $\BB_{av}$ and $\BB_{in}$ whose critical points
correspond to the Levi-Civit\`a regularizations of frozen planet orbits for the mean
and instantaneous interaction, respectively.
We proved that for each $r\in[0,1]$ the $L^2$-gradient $\nabla\BB_r$
of the interpolation 
$\BB_r=r\BB_{in}+(1-r)\BB_{av}$ is a
$C^1$-Fredholm map of index $0$. Now to each $C^1$-Fredholm map $F$ of
index $0$ with compact zero set one can associate its {\em mod $2$ Euler number}
$\chi(F)\in\Z/2\Z$ counting its zeroes modulo $2$ after perturbation,
see~\cite[Appendix C]{cieliebak-frauenfelder-volkov}. 
Using the compactness result from~\cite{frauenfelder2} and homotopy
invariance of the mod $2$ Euler number, we deduced
\begin{equation}\label{eq:chi-mod2}
   \chi(\nabla\BB_{in}) \equiv \chi(\nabla\BB_{av}) \equiv
   1\quad\text{(mod $2$)}.
\end{equation}
Here the functionals are considered on a suitable space of
normalized simple symmetric loops. 
In particular, for each negative energy there exists a frozen planet
orbit~\cite[Corollary C]{cieliebak-frauenfelder-volkov}. 

The proof that $\chi(\nabla\BB_{av})\equiv 1$ mod $2$ in~\cite[Theorem
  D.1]{cieliebak-frauenfelder-volkov} was based on a further
deformation of $\nabla\BB_{av}$ through Fredholm maps which are not
gradients of functionals. Our first result in the present paper
improves this to (see Corollary~\ref{cor:uniqueness-mean})

\textbf{Theorem\,A: }
{\it For each negative energy
there exists a unique normalized simple symmetric frozen planet orbit
for the mean intersection functional 
$\BB_{av}$. It is nondegenerate and of Morse index $0$.}
\smallskip

To prove this, we introduce a family of functionals $\FF_r$,
$r\in[0,\infty)$, such that for 
$\rho=(\sqrt{2}-1)^2$ the critical
points of $\FF_\rho$ and their Hessians agree with those of $\BB_{av}$, see 
Section~\ref{ss:modkeylemmas}.
On the other hand, $\FF_0$ describes the {\em free fall} of an
electron into the helium nucleus (again undergoing consecutive
collisions which are regularized), which is easily seen to possess a
unique simple periodic orbit that is nondegenerate of Morse index $0$.
This can be seen by looking at its Fourier expansion, 
see~\cite[Lemma~3.6]{frauenfelder-weber}.
Thus Theorem A (as well as the uniqueness of the frozen
planet orbit for $\BB_{av}$) will be a consequence of the following
result (see Theorem~\ref{nondeg}):
\smallskip

\textbf{Theorem\,B: }
{\it For every $r \in [0,\infty)$, each critical point of
    $\mathcal{F}_r$ is nondegenerate.} 
\smallskip

In Section~\ref{sec:uniqueness}
we explain how requiring the additional
properties ``normalized'', ``simple'', ``symmetric'' breaks the symmetries and leads to uniqueness.
The main ingredient in the proof of Theorem B is an algebraic identity
associated to critical points of $\FF_r$ which can be solved in terms
of elliptic integrals (see Proposition~\ref{prop:eq} and
Appendix~\ref{sec:ell-int}). 

Theorem A allows us to upgrade equation~\eqref{eq:chi-mod2} to an
equality of {\em integer valued Euler numbers}. For this, we need to
define a $\Z$-valued Euler number for a class of Fredholm sections 
including the $\nabla\BB_r$ above. The important feature of these
Fredholm sections is that their linearizations are {\em self-adjoint}
and bounded from below with respect to the $L^2$-scalar product. 
The main ingredient is the following abstract result which may be of
independent interest.

Let $F$ be a real Hilbert space, and $E\subset F$ a dense linear
subspace which is itself a Hilbert space (with a different inner
product) such that the inclusion $E\into F$ is compact. 
For a real number $\mathfrak{R}$ consider the spaces
$$
   \FF_s^{>\mathfrak{R}}(E,F) \subset \FF_s(E,F) \subset \FF(E,F)
$$
where $\FF(E,F)$ is the space of Fredholm operators $E\to F$,
$\FF_s(E,F)$ the subspace of operators that are self-adjoint as
unbounded operators on $F$ with domain $E$, and
$\FF_s^{>\mathfrak{R}}(E,F)$ the subspace of operators whose spectrum
is contained in $(\mathfrak{R},\infty)$. 
Recall (see e.g.~\cite{salamon}) that the determinants
$
   \det(D) = \Lambda^{\rm max}(\ker D^*)\otimes\Lambda^{\rm max}(\ker D)
$
of $D\in \FF(X,Y)$ give rise to a real line bundle, the {\em determinant line bundle}
$
   \det \to \FF(X,Y). 
$
It is well-known that this line bundle is non-orientable, but we have
(see Theorem~\ref{thm:det}): 
\smallskip

\textbf{Theorem\,C: }
{\it The restriction of the determinant line bundle to
$\FF_s^{>\mathfrak{R}}(E,F)$ carries a canonical orientation.} 
\smallskip

Contrarily to our initial expectation, the restriction of the
determinant line bundle to $\FF_s(E,F)$ is non-orientable; we
construct an explicit loop over which the bundle is nontrivial in
Proposition~\ref{prop:counterex}. 

Consider now a Hilbert manifold $X$ and a Hilbert space bundle $E\to X$
with a continuous bundle inclusion $TX\subset E$ such that $T_xX\subset
E_x$ is dense and the inclusion $T_xX\into E_x$ is compact for each
$x\in X$. Denote by ${\rm Func}_s^{>\mathfrak{R}}(X)$ the space of
$C^2$-functions $f:X\to\R$ whose $E$-gradient $\Nabla{E}f:X\to E$ is of
class $C^1$ such that $\Nabla{E}f(x)\in\FF_s^{>\mathfrak{R}}(T_xX,E_x)$
for each critical point $x$ of $f$. For this class of functions,
Theorem C allows us to define a $\Z$-valued Euler number (see 
Theorem~\ref{thm:zero-Fredholm}): 
\smallskip

\textbf{Corollary\,A: }
To each $f\in{\rm Func}_s^{>\mathfrak{R}}(X)$ with compact zero set
we can associate an {\em Euler number} $\chi(\Nabla{E}f)\in\Z$ 
which is uniquely characterized by suitable axioms of (Transversality), 
(Excision), and (Homotopy).
\smallskip

Verifying that the $\BB_r$ above belong to
${\rm Func}_s^{>\mathfrak{R}}(X)$ for a suitable bundle $E\to X$ and
using the (Transversality) and (Homotopy) axioms, we deduce (see
Corollary~\ref{cor:frozen-planet-count}) 

\textbf{Corollary\,B: }
The integral count of normalized simple symmetric frozen planet orbits equals
\begin{equation*}%\label{eq:chi}
   \chi(\nabla\BB_{in})  = \chi(\nabla\BB_{av}) = 1 \in \Z.
\end{equation*}

\textbf{Remark: }
Symmetric frozen planet orbits have two Morse indices: one as
a symmetric frozen planet orbit, and one just as a frozen planet orbit
forgetting about the symmetry. When we talk about Morse indices (which
enter into the Euler number via the (Transversality) axiom in
Theorem~\ref{thm:zero-Fredholm}) we always mean the symmetric
one. In fact, the two indices are different: the unique normalized
simple symmetric periodic orbit for the regularized free fall (i.e., 
the functional $\mathcal{F}_r$ for $r=0$) has index $0$ as a symmetric
orbit and index $1$ just as a periodic orbit. As we show, the functional
$\mathcal{F}_r$ is always nondegenerate in the symmetric as well as in
the just periodic sense. Therefore, the symmetric Euler characteristic
equals $1$ while the periodic Euler characteristic equals
$-1$. However, to our knowledge, compactness for frozen planet orbits
in the homotopy from mean to instantaneous interaction 
has only been established in the symmetric case~\cite{frauenfelder2}.
Therefore, it is not clear whether Corollary\,B has an analogue if
one forgets about the symmetry. We expect that it does and that the
Euler characteristic for normalized simple (not necessarily symmetric)
frozen planet orbits with instantaneous interaction is $-1$. This
would fit with the findings of physicists \cite{tanner-richter-rost,
  wintgen-richter-tanner}, who numerically detected a frozen planet
orbit for instantaneous interaction which is stable and therefore has
odd Conley-Zehnder index. 

{\bf Acknowledgements. }
%marginpar{\tiny Is this o.k.? Anyone else?}
We thank Bernd Schmidt for the elegant proof of Lemma~\ref{lem:Bernd}, and
Dirk Bl\"omker for helpful conversations on Dunford calculus. 
This research was supported by Deutsche Forschungsgemeinschaft (DFG,
German Research Foundation) Projects CI 45/8-2 and FR 2637/2-2.

%%%%%%%%%%%%%%%%%%%%%%%%%%%%%%%%%%%%%%%%%%%%%%%%%%%%%%%%%%%%%%%%
\section{Levi-Civita transformation}\label{sec:LC}
%%%%%%%%%%%%%%%%%%%%%%%%%%%%%%%%%%%%%%%%%%%%%%%%%%%%%%%%%%%%%%%%

In this section we recall some background on the Levi-Civita
transformation. For details we refer the reader 
to~\cite{cieliebak-frauenfelder-volkov}.
%idea from~\cite{barutello-ortega-verzini}
%without worrying about the regularity of the involved maps; precise
%statements will be given in the following subsections.

We abbreviate by $S^1=\mathbb{R}/\mathbb{Z}$ the circle. We denote the
$L^2$-inner product of $z_1,z_2\in L^2(S^1,\mathbb{R})$ by
$$
  \langle z_1,z_2\rangle := \int_0^1 z_1(\tau) z_2(\tau) d\tau,
$$
and the $L^2$-norm of $z \in L^2(S^1, \mathbb{R})$ by
$$
  \|z\| := \sqrt{\langle z, z\rangle}.
$$
In the sequel we will work with Sobolev spaces $H^k=W^{k,2}$, but
the only relevant norms and inner products will be the ones from $L^2$.

Consider two maps
$$
   q:S^1\to\R_{\geq 0},\qquad z:S^1\to\R
$$
related by the {\em Levi-Civita transformation}
\begin{equation}\label{eq:LC}
   q(t) = z(\tau)^2
\end{equation}
for a time change $t\longleftrightarrow \tau$ satisfying
$0\longleftrightarrow 0$ and
\begin{equation}\label{eq:t-tau}
   \frac{dt}{q(t)} = \frac{d\tau}{\|z\|^2}. 
\end{equation}
This implies that the mean values of $q$ and $1/q$ are given by
\begin{equation}\label{bov1}
  \overline{q} :=
  \int_0^1 q(t)dt
  = \int_0^1 \frac{z(\tau)^4}{\|z\|^2} d\tau
  = \frac{\|z^2\|^2}{\|z\|^2}
\end{equation}
and
\begin{equation}\label{bov2}
   \int_0^1\frac{dt}{q(t)} = \frac{1}{\|z\|^2}.
\end{equation}
We will denote derivatives with respect to $t$ by a dot and
derivatives with respect to $\tau$ by a prime. Then the first and
second derivatives of $q$ and $z$ (where they are defined) are related by
\begin{equation}\label{bov3}
  \dot q(t)
  = 2z(\tau)z'(\tau)\frac{d\tau}{dt} = \frac{2\|z\|^2z'(\tau)}{z(\tau)}
\end{equation}
and
\begin{equation}\label{eq:ddot-q-z}
  \ddot q(t)
  = 2\|z\|^2\frac{z''(\tau)z(\tau)-z'(\tau)^2}{z(\tau)^2}\frac{d\tau}{dt}
  = \frac{2\|z\|^4}{z(\tau)^4}\bigl(z''(\tau)z(\tau)-z'(\tau)^2\bigr).
\end{equation}
Substituting $z^2$ and $z'^2$ by~\eqref{eq:LC} and~\eqref{bov3} this becomes
\begin{equation}\label{eq:ddot-q}
  \ddot q(t)
  = \frac{1}{q(t)}\Bigl(2\|z\|^4\frac{z''(\tau)}{z(\tau)} - \frac{\dot q(t)^2}{2}\Bigr).
\end{equation}
The $L^2$-norm of the derivative of $q$ is given by 
\begin{eqnarray}\label{bov4}
  \|\dot{q}\|^2
  = \int_0^1 \dot{q}(t)^2 dt %\\ \nonumber
  = \int_0^1
  \frac{4\|z\|^4z'(\tau)^2}{z(\tau)^2}\,\frac{z(\tau)^2}{\|z\|^2}d\tau
  %\\ \nonumber
  = 4\|z\|^2\|z'\|^2.
\end{eqnarray}

We can now give the precise definition of the Levi-Civita transformation.
Let $z\in C^0(S^1,\R)$ be a continuous function with finite zero set
$$
  Z_z:=z^{-1}(0). 
$$
We associate to $z$ a $C^1$-map $t_z:S^1\to S^1$ by
\begin{equation}\label{eq:tz}
   t_z(\tau) := \frac{1}{\|z\|^2}\int_0^\tau z(\sigma)^2d\sigma.
\end{equation}
Note that $t_z(0)=0$ and
\begin{equation}\label{eq:dert}
   t_z'(\tau) = \frac{z(\tau)^2}{\|z\|^2}.
\end{equation}
Since $z$ has only finitely many zeroes, this shows that $t_z$ is
strictly increasing and we conclude

\begin{lemma}\label{lem:tz-homeo}
If $z\in C^0(S^1,\R)$ has only finitely many zeroes, then
the map $t_z:S^1\to S^1$ defined by~\eqref{eq:tz} is a homeomorphism. 
\hfill$\square$
\end{lemma}

It follows that $t_z \colon S^1 \to S^1$ has a continuous inverse
$$
  \tau_z:=t_z^{-1} \colon S^1 \to S^1.
$$
Since $t_z$ is of class $C^1$, the function $\tau_z$ is also of class $C^1$ on the
complement of the finite set $t_z(Z_z)$ with derivative
\begin{equation}\label{eq:taudot}
  \dot\tau_z(t) = \frac{\|z\|^2}{z(\tau_z(t))^2},\qquad t\in
  S^1\setminus t_z(Z_z).
\end{equation}
We define a continuous map $q:S^1\to \R_{\ge 0}$ by
\begin{equation}\label{eq:LC-def}
  q(t) := z(\tau_z(t))^2.
\end{equation}
Then the two maps $z,q$ are related by the Levi-Civita
transformation~\eqref{eq:LC} with $\tau=\tau_z$. Their zero sets 
$$
  Z_z=z^{-1}(0)\quad \text{and}\quad Z_q:=q^{-1}(0)=t_z(Z_z)
$$ 
are in bijective correspondence via $t_z$ (or equivalently $\tau_z$).
Moreover, by~\eqref{bov2} we have
$$
   \int_0^1\frac{ds}{q(s)} = \frac{1}{\|z\|^2} <\infty.
$$
Conversely, suppose we are given a map $q\in C^0(S^1,\R_{\geq 0})$
with finite zero set $Z_q$ satisfying $\int_0^1\frac{ds}{q(s)}<\infty$. 
We associate to $q$ the time reparametrization $\tau_q:S^1\to S^1$,
\begin{equation}\label{eq:tauq}
   \tau_q(t) :=
   \Bigl(\int_0^1\frac{ds}{q(s)}\Bigr)^{-1} \int_0^t\frac{1}{q(s)}ds. 
\end{equation}
Then $\tau_q(1)=1$, $\tau_q$ is of class $C^1$ outside the zero set $Z_q=q^{-1}(0)$
with derivative
\begin{equation}\label{eq:tauqdot}
   \tau_q'(t) =
   \Bigl(\int_0^1\frac{ds}{q(s)}\Bigr)^{-1}\frac{1}{q(t)}, \qquad t\in
   S^1\setminus Z_q.
\end{equation}
By~\cite[Lemma 2.1]{barutello-ortega-verzini}, the map $\tau_q:S^1\to S^1$
is a homeomorphism whose inverse $t_q:=\tau_q^{-1}$ is of class $C^1$
and satisfies $t_q(1)=\tau_q^{-1}(1)=1$ and
\begin{equation}\label{eq:tqdot}
   t_q'(\tau) = \Bigl(\int_0^1\frac{ds}{q(s)}\Bigr)\, q(t_q(\tau)),\qquad
   \tau\in S^1.
\end{equation}
Suppose that $z:S^1\to\R$ is a continuous function satisfying
\begin{equation}\label{eq:LC-inv}
   z(\tau)^2 = q(t_q(\tau)).
\end{equation}
Then $z$ has finite zero set $Z_z=\tau_q(Z_q)$, so we can associate to
$z$ the homeomorphism $t_z:S^1\to S^1$ defined by~\eqref{eq:tz} and
its inverse $\tau_z$. We claim that 
\begin{equation}\label{eq:homeoqz}
   \tau_q=\tau_z\quad \text{and}\quad t_q=t_z.
\end{equation}
It is enough to check the second equality. For this we compute 
$$
  \int_0^\tau z(\sigma)^2d\sigma
  =\int_0^\tau q(t_q(\sigma))d\sigma
  \stackrel{(*)}{=}\Bigl(\int_0^1\frac{ds}{q(s)}\Bigr)^{-1}
  \int_0^{t_q(\tau)}ds
  =\Bigl(\int_0^1\frac{ds}{q(s)}\Bigr)^{-1}t_q(\tau),
$$
where $(*)$ follows from the coordinate change $\sigma=\tau_q(s)$
and~\eqref{eq:tauqdot}. Evaluating at $\tau=1$ gives us
\begin{equation}\label{eq:meanqinv}
   \frac{1}{\|z\|^2}=\int_0^1\frac{ds}{q(s)}.
\end{equation}
Therefore,
$$
   t_z(\tau) = \frac{1}{||z||^2}\int_0^\tau z(\sigma)^2d\sigma = t_q(\tau)
$$
and~\eqref{eq:homeoqz} is established. Hence $q$ is the Levi-Civita
transform of $z$ defined by~\eqref{eq:LC-def}. 

Equation~\eqref{eq:LC-inv} does not uniquely determine $z$ for given
$q$ because the sign of $z$ can be arbitrarily chosen on each connected
component of $S^1\setminus Z_z$.
If $Z_z$ consists of an {\em even} number of points, then we can
determine $z$ up to a global sign by the requirement that $z$ changes
its sign at each zero. 
If $Z_z$ consists of an {\em odd} number of points, then the
requirement that $z$ changes its sign at each zero leads to
$z(\tau+1)=-z(\tau)$, so $z$ has period $2$ rather than $1$. 
Therefore, the preceding discussion shows

\begin{lemma}\label{lem:LC}
The Levi-Civita transformation $z\mapsto q$ given by~\eqref{eq:LC-def}
defines for each {\em even} integer $m$ a surjective 2-to-1 map
\begin{align*}
   \LL\colon &\{z\in C^0(S^1,\R)\mid \text{$z$ has precisely $m$ zeroes and
     switches sign at each zero}\} \\
   &\longrightarrow \{q\in C^0(S^1,\R_{\geq 0})\mid \text{$z$ has precisely $m$
     zeroes and }\int_0^1\frac{ds}{q(s)}<\infty\},
\end{align*}
and for each {\em odd} integer $m$ a surjective 2-to-1 map
\begin{align*}
   \LL\colon &\{z\in C^0(\R/2\Z,\R)\mid \text{$z$ has precisely $2m$ zeroes and
     switches sign at each zero},\\
   &\ \ \ \ \ \ \ \ \ \ \ \ \ \ \ \ \ \ \ \ \ \ \ \ z(\tau+1)=-z(\tau)
   \text{ for all }\tau\} \\
   &\longrightarrow \{q\in C^0(S^1,\R_{\geq 0})\mid \text{$z$ has precisely $m$
     zeroes and }\int_0^1\frac{ds}{q(s)}<\infty\}.
\end{align*}
\hfill$\square$
\end{lemma}

\section{The functionals $\FF_r$ and their critical points}\label{sec:frozen}
%%%%%%%%%%%%%%%%%%%%%%%%%%%%%%%%%%%%%%%%%%%%%%%%%%%%%%%%%%%%%%%%

We denote by $S^1=\mathbb{R}/\mathbb{Z}$ the circle and 
abbreviate by 
$$
  H^1_*(S^1,\mathbb{R})=H^1(S^1,\mathbb{R})\setminus \{0\},\quad
  H^2_*(S^1,\mathbb{R})=H^2(S^1,\mathbb{R})\setminus \{0\} 
$$ 
the open subsets of the Hilbert space $H^1(S^1,\mathbb{R})$ 
respectively $H^2(S^1,\mathbb{R})$ where the origin is removed. For $r
\in [0,\infty)$ we consider the functional
$$
  \mathcal{F}_r \colon H^1_*(S^1,\mathbb{R}) \to \mathbb{R}, \quad 
  z \mapsto 2||z||^2||z'||^2+\frac{2}{||z||^2}+r\frac{||z||^2}{||z^2||^2},
$$
where as before $||z||$ is the $L^2$-norm of the loop $z$. If $z \in
H^1_*(S^1,\mathbb{R})$ and $\xi \in H^1(S^1,\mathbb{R})$, then the
differential of $\mathcal{F}_r$ at $z$ in direction of $\xi$ is given by
\begin{eqnarray*}
D\mathcal{F}_r(z)\xi&=&4||z||^2\langle z',\xi'\rangle+4||z'||^2 \langle z,\xi\rangle
-\frac{4\langle z,\xi\rangle}{||z||^4}+r
\bigg(\frac{2\langle z,\xi\rangle}{||z^2||^2}-\frac{4 ||z||^2\langle z^3,\xi\rangle}{||z^2||^4}\bigg),
\end{eqnarray*}
where $\langle\,,\,\rangle$ denotes the $L^2$-inner product.
Therefore, integration by parts shows that for 
$z\in H^2_*(S^1,\mathbb{R})$ the functional 
$\mathcal{F}_r$ possess an $L^2$-gradient given by the 
formula
\begin{eqnarray}\label{eq:grad}
\nabla\mathcal{F}_r(z)&=&-4\|z\|^2(z''+bz+2az^3),\\ \nonumber
b&=&\frac{1}{||z||^6}-\frac{||z'||^2}{||z||^2}-\frac{r}{2||z||^2\cdot ||z^2||^2},\\ \nonumber
a&=&\frac{r}{2||z^2||^4}.
\end{eqnarray}
We see that $\nabla\mathcal{F}_r:H^2_*(S^1,\R)\to L^2(S^1,\R)$ is
differentiable; we call its derivative the {\em Hessian} of $\mathcal{F}_r$.
Critical points of $\mathcal{F}_r$ are solutions of the problem
\begin{equation}\label{eq:crit}
z''=-bz-2az^3.
\end{equation}
A standard bootstrapping argument implies
that critical points of $\mathcal{F}_r$ are in fact 
smooth. If $z$ is a solution of (\ref{eq:crit}), then for $n \in \mathbb{N}$ the loop $z_n$ defined as
$$z_n(\tau)=n^{-1/3}z(n\tau), \quad \tau \in S^1$$
is another solution of (\ref{eq:crit}) with $a_n=n^{8/3}a$ and $b_n=n^2 b$. We say that a critical point
$z$ is \emph{multiply covered} if there exists a critical point $w$ and $n>1$ such that $z=w_n$.
Otherwise we call the critical point \emph{simple}. Moreover, the functional $\mathcal{F}_r$ is invariant under the $S^1$-action given by time shift, 
$$\sigma_*z(\tau)=z(\tau+\sigma), \quad \tau \in S^1,$$
where $\sigma \in S^1$ and $z \in H^1_*(S^1,\mathbb{R})$. In
particular, its critical points are invariant under time shift as well. Therefore, if $z$ is a critical point of $\mathcal{F}_r$,
then $z'$ lies in the kernel of the Hessian of $\mathcal{F}_r$ at $z$. 
\\ \\
%We say that a critical point $z$ of $\mathcal{F}_r$ is of
%\emph{collision type} if
It follows from Proposition~\ref{prop:z-to-q} below that for each critical
point $z$ of $\mathcal{F}_r$, $r\geq 0$, there exists
$\tau_0 \in S^1$ such that $z(\tau_0)=0$. 
Indeed, otherwise there exists 
$t_0\in S^1$, such that the Levi-Civita 
transform $q$ of $z$ attains a minimum
$q(t_0)>0$ at $t_0$. Therefore,
$\ddot{q}(t_0)\ge 0$. On the other hand, according to 
Proposition~\ref{prop:z-to-q}, the function $q$ satisfies 
equation~\eqref{eq:diffeq-q}. 
For $t=t_0$
this implies that $\ddot{q}(t_0)<0$.
We then necessarily have
$z'(\tau_0) \neq 0$, since otherwise by (\ref{eq:crit}) the loop $z$
would be the constant loop at the origin which does not lie in 
$H^1_*(S^1,\mathbb{R})$. Thus a critical point of $\FF_r$ is never a
fixed point of the $S^1$-action on the free loop space. In particular,
its nullity, i.e., the dimension of the kernel of its Hessian, is at
least one. We say that a critical point is \emph{nondegenerate} if its
nullity is precisely one.
%It follows from Theorem~\ref{thm:z-to-q} below that for nonnegative
%$r$, all critical points are of collision type.
Our first result asserts that nondegeneracy always holds true for
nonnegative $r$. It corresponds to Theorem B from the Introduction
and will be proved in Section~\ref{sec:nondeg}.

\begin{thm}\label{nondeg}
For every $r \in [0,\infty)$, each critical point of $\mathcal{F}_r$ is nondegenerate.
\end{thm}

%If $z$ is a critical point of collision type, then its iterates $z_n$
%are critical points of collision type as well. In particular, the
%theorem tells us that for a critical point of collision type all its
%iterates are nondegenerate.  

%%%%%%%%%%%%%%%%%%%%%%%%%%%%%%%%%%%%%%%%%%%%%%%%%%%%%%%%%%%%%%%%
\subsection{Levi-Civita transform of critical points}\label{sec:qpic}
%%%%%%%%%%%%%%%%%%%%%%%%%%%%%%%%%%%%%%%%%%%%%%%%%%%%%%%%%%%%%%%%

In this section we apply the Levi-Civita transformation to critical
points of $\mathcal{F}_r$. Let $z \in H^1_*(S^1,\mathbb{R})$ be a
solution of~\eqref{eq:crit}
%of collision type, i.e., a solution of~\eqref{eq:crit} for which there exists $\tau_0 \in S^1$ such that $z(\tau_0)=0$.
and $q(t):=z(\tau_z(t))^2$ its Levi-Civita transform. 
We compute at points $t\in S^1\setminus t_{z}(Z_z)$:
\begin{eqnarray}
\ddot{q} \nonumber
&\stackrel{(A)}{=}&\frac{1}{q}\left(-2||z^4||(b+2az^2)-
\frac{\dot q^2}{2}\right)\\ \nonumber
&\stackrel{(B)}{=}  &\frac{1}{q}
\left(-\frac{2}{||z||^2}+2||z||^2||z'||^2+
r\frac{||z||^2}{||z^2||^2}-2r\frac{||z||^4}
{||z^2||^4}q-\frac{\dot q^2}{2}\right)
\\ \nonumber
&\stackrel{(C)}{=}& \frac{1}{q}\left(-\int_0^1\frac{2}{q(s)}ds+
\frac{||\dot q||^2}{2}+\frac{r}{\overline{q}}-2r\frac{q}{\overline{q}^2}-\frac{\dot q^2}{2}\right).
\end{eqnarray}
Equality~(A) follows from substituting $z''$ by~\eqref{eq:crit} in equation~\eqref{eq:ddot-q};
equality~(B) uses the expressions for $a$ and $b$
in~\eqref{eq:crit}; equality~(C) uses $q=z^2$ as well
as~\eqref{eq:t-tau},~\eqref{bov1} and~\eqref{bov4}.
Thus $q$ satisfies the ODE
\begin{equation}\label{bov6}
  \ddot{q} = \bigg(c
  -\frac{\dot{q}^2}{2}\bigg)\frac{1}{q}-\frac{2r}{\overline{q}^2}
\end{equation}
with the constant
\begin{equation}\label{eq:c1-mean}
  c = \frac{||\dot{q}||^2}{2}-\int_0^1
  \frac{2}{q(s)}ds+\frac{r}{\overline{q}}. 
\end{equation}
At the global maximum $t_{\rm max}$ of $q$, 
equation~\eqref{bov6} becomes 
$$
  \frac{c}{q(t_{\rm max})} + 
  \frac{2r}{\overline{q}^2} =
  \ddot q_1(t_{\rm max})\leq 0,
$$
hence
\begin{equation}\label{eq:c1-ineq}
   c \leq -\frac{2rq(t_{\rm max})}
   {\overline{q}^2}.
\end{equation}
Let now $t_-<t_+$ be adjacent zeroes of $q$ and consider the smooth map
$$
  \beta:=\frac{\ddot{q}+\frac{r}{\overline{q}^2}}{q}
  \colon (t_-,t_+) \to \mathbb{R}.
$$
From (\ref{bov6}) we obtain
$$
  \beta q^2
  = c -\frac{\dot{q}^2}{2} -
  \frac{rq}{\overline{q}^2}.
$$
With inequality~\eqref{eq:c1-ineq} this implies
$$
  \beta q^2
  \leq -\frac{2rq(t_{\rm max})}
   {\overline{q}^2}  -\frac{\dot{q}^2}{2} -
  \frac{rq}{\overline{q}^2}  < 0,
$$
hence $\beta<0$ on $(t_-,t_+)$. 
Differentiating both sides of the
equation for $\beta q^2$ we get
$$\dot{\beta} q^2+2 \beta q \dot{q}=-\ddot{q}\dot{q}
-\frac{r\dot q}{\overline{q}^2}
=-\beta q \dot{q},$$
and therefore
\begin{equation}\label{eq:betaqu}
  \dot{\beta} q = -3\beta \dot{q}.
\end{equation}
We need the following
\begin{lemma}\label{lem:betamu}
Equation~\eqref{eq:betaqu} for functions $q>0$ and 
$\beta<0$ on $(t_-,t_+)$ implies that
\begin{equation}\label{eq:betamu}
  \beta = -\frac{\mu}{q^3}
\end{equation}
on $(t_-,t_+)$ for some constant $\mu>0$.
\end{lemma}

\textbf{Proof: }
Dividing both sides of equation~\eqref{eq:betaqu} by $q\beta$ yields
$$
  \ddt \log(-\beta)=-3\ddt \log(q),
$$
which by integration implies the lemma.
\hfill $\square$

By Lemma~\ref{lem:betamu} we get 
equation~\eqref{eq:betamu}
on $(t_-,t_+)$ for some constant $\mu>0$.
By definition of $\beta$, this yields the following equation for $q$:
\begin{equation}\label{bov7}
  \ddot{q}(t) =
  -\frac{\mu}{q(t)^2}-\frac{r}{\overline{q}^2}
\end{equation}
for $t\in (t_-,t_+)$. 
It remains to compute $\mu$.
Plugging this into~\eqref{bov6} we infer
\begin{equation}\label{bov7a}
\mu= -\bigg(c-\frac{\dot{q}(t)^2}{2}\bigg)q(t)
+\frac{rq(t)^2}{\overline{q}^2}
\end{equation}
for $t\in(t_-,t_+)$. In particular, using~\eqref{bov3} and
$q(t_\pm)=0$ we obtain 
\begin{equation}\label{bov8}
  \mu = \lim_{t \to t_\pm} \frac{\dot{q}(t)^2 q(t)}{2}
  = 2||z||^4 z'\big(\tau_{z}(t_\pm)\big)^2.
\end{equation}
We deduce from this that equation~\eqref{bov7} holds on $S^1 \setminus
t_z(Z_z)$ with a fixed $\mu$ independent of 
the connected component in $S^1 \setminus
t_z(Z_z)$. Dividing~\eqref{bov7a} by $q(t)$ and inserting $c$
from~\eqref{eq:c1-mean}, we get
\begin{equation}\label{bov7b}
  \frac{\mu}{q(t)}
  = -\frac{\|\dot{q}\|^2}{2}+\int_0^1
  \frac{2}{q(s)}ds-\frac{r}{\overline{q}} 
  +\frac{\dot{q}(t)^2}{2}+
  \frac{rq(t)}{\overline{q}^2}.
\end{equation}
Integrating this equation yields
$$
  \mu\int_0^1 \frac{1}{q(t)}dt=2\int_0^1 \frac{1}{q(s)}ds,
$$
and therefore
$$
\mu=2.
$$
This gives us the following statement.

\begin{prop}\label{prop:z-to-q}
Assume that $z\in H^1_*(S^1,\mathbb{R})$
is a critical point of the frozen functional 
$\mathcal{F}_r$. Then the Levi-Civita transform
$q(t)=z(\tau_z(t))^2$ of $z$ satisfies the 
differential equation
\begin{equation}\label{eq:diffeq-q}
  \ddot{q}(t) =
  -\frac{2}{q(t)^2}-\frac{r}{\overline{q}^2}.
\end{equation}
\end{prop}

Equation~\eqref{eq:diffeq-q} explains the physical meaning of critical
points of $\FF_r$: The orbit $q(t)$ describes an electron on the line
attracted by a doubly positively charged nucleus at the origin and
subject to an additional force depending on its mean $\ol q$.
For $r=0$ the mean interaction force vanishes and $q(t)$ describes
the {\em free fall} of the electron into the nucleus, regularized by
elastic reflection as it hits the nucleus.
For $r>0$ the mean interaction force pushes the electron towards the
origin, which can be interpreted as the effect of a second electron further
away and on the same side of the nucleus. Indeed, we will show in
\S\ref{sec:appl-frozen} that for a suitable value of $r$ the system
describes frozen planet orbits in helium with mean interaction between
the electrons.
For $r<0$ the mean interaction force pushes the electron away from the
origin, which can be interpreted as the effect of a second electron 
on the other side of the nucleus. However, we will not consider the
case $r<0$ in this paper.

%%%%%%%%%%%%%%%%%%%%%%%%%%%%%%%%%%%%%%%%%%%%%%%%%%%%%%%%%%%%%%%%
\subsection{Analysis of critical points}
%%%%%%%%%%%%%%%%%%%%%%%%%%%%%%%%%%%%%%%%%%%%%%%%%%%%%%%%%%%%%%%%

Let $r\geq 0$ and $z \in H^1_*(S^1,\mathbb{R})$ be a critical point of
$\mathcal{F}_r$,
%of collision type,
i.e., a solution of~\eqref{eq:crit}. 
%marginpar{Suitable language like ``time reparametrization''
%here to be added where needed. EV.}
We denote by $||z||_0$ the maximum norm of $z$ and introduce the quantities
\begin{equation}\label{eq:vw}
v:=\frac{||z||^2}{||z||_0^2}, \qquad w:=\frac{||z||^2 ||z||_0^2}{||z^2||^2}.
\end{equation}
For $n \in \mathbb{N}_0$ we consider the elliptic integral
$$I_n \colon (-\infty,1) \to \mathbb{R}.\qquad I_n(m):=\int_0^1
\frac{\zeta^{2n}}{\sqrt{(1-\zeta^2)(1-m\zeta^2)}}d\zeta.$$ 
The main result of this section is the following proposition.

\begin{prop}\label{prop:eq}
The quantities $v$ and $w$ satisfy the equations
\begin{equation}\label{eq}
v=\frac{2}{4+3rw-2rw^2}=\frac{I_1}{I_0}\bigg(-\frac{rw^2}{2}\bigg).
\end{equation}
\end{prop}

\textbf{Proof: }
After a time reparametrization, we can assume without loss of generality
that $z$ attains at time $0$ its global maximum
$$z_0:=||z||_0=z(0).$$
Recall that there exists $\tau_0
\in S^1$ such that $z(\tau_0)=0$. 
We further choose $\tau_0$ as the smallest number in $(0,1)$ with the property
that $z(\tau_0)=0$.

We first eliminate in the formula for $b$ in~\eqref{eq:grad} the variable $||z'||$. 
Using~\eqref{eq:crit}, we obtain via integration by parts
$$||z'||^2=-\langle z,z''\rangle=b||z||^2+2a||z^2||^2.$$
Plugging this into the equation for $b$ and using the equation for
$a$, we get
$$b=\frac{1}{||z||^6}-b-\frac{r}{||z||^2 \cdot ||z^2||^2}-\frac{r}{2||z||^2\cdot ||z^2||^2},$$
implying 
\begin{equation}\label{beq}
b=\frac{1}{2||z||^6}-\frac{3r}{4||z||^2\cdot ||z^2||^2}.
\end{equation}
From~\eqref{eq:crit} we further infer that we have the conserved quantity
\begin{equation}\label{pres}
\frac{z'(\tau)^2}{2}+\frac{bz(\tau)^2}{2}+\frac{az(\tau)^4}{2}=\frac{c}{2}
\end{equation}
for some constant $c$. Since $z(\tau_0)=0$ and $z'(\tau_0) \neq 0$ we conclude that $c$ is positive. 
Since $z$ attains its global maximum $z_0$ at time $\tau=0$, we necessarily have $z'(0)=0$ and therefore
$$b z_0^2+a z_0^4=c.$$
This means that $z_0^2$ is a root of the quadratic polynomial
$$p(x)=ax^2+bx-c.$$
Moreover, positivity of $c$ yields
\begin{equation}\label{eq:c-pos}
  az_0^2+b>0.
\end{equation} 
\textbf{Case\,1:} $r > 0$.
\\ \\
In this case $a \neq 0$ as well and the second root of
$p$ is given by $-z_0^2-\tfrac{b}{a}$. In particular, the quadratic polynomial factorizes as
$$p(x)=\Big(x-z_0^2\Big)\Big(ax+az_0^2+b\Big).$$
Plugging this into~\eqref{pres} we obtain
\begin{equation}\label{poly}
z'(t)^2=-p\bigl(z(\tau)\bigr)=\Big(z_0^2-z(\tau)^2\Big)\Big(a z(\tau)^2+az_0^2+b\Big), \quad \tau \in S^1.
\end{equation}
Equation~\eqref{eq:crit} is invariant under reflection at the origin and time reversal. 
Since $z(\tau_0)=0$, we conclude that
$$z(\tau)=-z(2\tau_0-\tau), \quad \tau \in [\tau_0,2\tau_0].$$
In particular, we have $z(2\tau_0)=-z(0)$ and $z'(2\tau_0)=0$. Using once more invariance under time reversal, we conclude that
$$z(\tau)=-z(2\tau_0+\tau), \quad \tau \in [0,4\tau_0].$$
In particular, we have $z(4\tau_0)=z(0)$ and $z'(4\tau_0)=z'(0)=0$. We see that $z$ is periodic
of period $4\tau_0$. Since $\tau_0$ was chosen as the first positive time at which $z$ passes through the
origin, we conclude that $4\tau_0$ is the minimal period of $z$. Since
$z$ is by assumption periodic of period
$1$, we conclude that there exists $n \in \mathbb{N}$ such that $4\tau_0
n=1$, i.e.,
$$t_0=\frac{1}{4n}.$$
If $z$ is simple, then $n=1$; otherwise it is multiply covered. Using (\ref{poly})
we therefore obtain
\begin{eqnarray}\label{e0}
\frac{1}{4n}&=&\int_0^{z_0}\frac{1}{\sqrt{(z_0^2-z^2)(az^2+az_0^2+b)}}dz\\ \nonumber
&=&\frac{1}{\sqrt{az_0^2+b}}\int_0^1\frac{1}{\sqrt{1-\zeta^2)\big(1+\frac{az_0^2}{az_0^2+b}\zeta^2\big)}}d\zeta\\ \nonumber
&=&\frac{1}{\sqrt{az_0^2+b}}I_0\bigg(-\frac{az_0^2}{az_0^2+b}\bigg).
\end{eqnarray}
Similarly, we compute
\begin{eqnarray}\label{e1}
\frac{||z||^2}{4n}&=&\int_0^{z_0}\frac{z^2}{\sqrt{(z_0^2-z^2)(az^2+az_0^2+b)}}dz\\ \nonumber
&=&\frac{z_0^2}{\sqrt{az_0^2+b}}\int_0^1\frac{\zeta^2}{\sqrt{1-\zeta^2)\big(1+\frac{az_0^2}{az_0^2+b}\zeta^2\big)}}d\zeta\\ \nonumber
&=&\frac{z_0^2}{\sqrt{az_0^2+b}}I_1\bigg(-\frac{az_0^2}{az_0^2+b}\bigg)
\end{eqnarray}
and
\begin{eqnarray}\label{e2}
\frac{||z^2||^2}{4n}&=&\int_0^{z_0}\frac{z^4}{\sqrt{(z_0^2-z^2)(az^2+az_0^2+b)}}\\ \nonumber
&=&\frac{z_0^4}{\sqrt{az_0^2+b}}\int_0^1\frac{\zeta^4}{\sqrt{1-\zeta^2)\big(1+\frac{az_0^2}{az_0^2+b}\zeta^2\big)}}d\zeta\\ \nonumber
&=&\frac{z_0^4}{\sqrt{az_0^2+b}}I_2\bigg(-\frac{az_0^2}{az_0^2+b}\bigg). 
\end{eqnarray}
The elliptic function $I_2$ can be expressed using the elliptic functions $I_0$ and $I_1$
by
$$I_2(m)=\frac{2(m+1)I_1(m)}{3m}-\frac{I_0(m)}{3m}$$
as explained in (\ref{i2}) in  Appendix~\ref{sec:ell-int}. Hence from (\ref{e0}), (\ref{e1}), and (\ref{e2}) we
obtain the equality
\begin{eqnarray*}
||z^2||^2&=&\frac{4nz_0^4}{\sqrt{az_0^2+b}}\Bigg(\frac{2\Big(\frac{az_0^2}{az_0^2+b}-1\Big)}{3\frac{az_0^2}{az_0^2+b}}I_1\bigg(-\frac{az_0^2}{az_0^2+b}\bigg)+\frac{az_0^2+b}{3az_0^2}I_0\bigg(-\frac{az_0^2}{az_0^2+b}\bigg)\Bigg)\\ 
&=&\frac{z_0^4}{\sqrt{az_0^2+b}}\Bigg(\frac{az_0^2+b}{3az_0^2}\sqrt{az_0^2+b}-
\frac{2b}{3az_0^2}\cdot \frac{\sqrt{az_0^2+b}}{z_0^2}||z||^2\Bigg)\\ 
&=&\frac{(az_0^2+b)z_0^2}{3a}-\frac{2b ||z||^2}{3a}\\ 
&=&\frac{z_0^4}{3}+\frac{b}{3a}\big(z_0^2-2||z||^2\big)\\ 
&=&\frac{z_0^4}{3}+\bigg(\frac{||z^2||^4}{3r||z||^6}-\frac{||z^2||^2}{2||z||^2}\bigg)
\cdot \big(z_0^2-2||z||^2\big)\\
&=&\frac{z_0^4}{3}+\frac{||z^2||^4 z_0^2}{3r||z||^6}-\frac{||z^2||^2z_0^2}{2||z||^2}
-\frac{2||z^2||^4}{3r||z||^4}+||z^2||^2
\end{eqnarray*}
where in the second to last equality we have used the equation for $a$ from \eqref{eq:grad} and 
equation~\eqref{beq} for $b$. Removing $||z^2||^2$ on both sides and multiplying the remaining terms by
$\frac{6r||z||^4}{||z^2||^4}$ we obtain the equation
$$0=\frac{2r ||z||^4 z_0^4}{||z^2||^4}+\frac{2z_0^2}{||z||^2}-\frac{3r||z||^2z_0^2}{||z^2||^2}-4.$$
By definition of $v$ and $w$ we can rewrite this as
$$0=2r w^2+\frac{2}{v}-3rw-4$$
or equivalently
$$v=\frac{2}{4+3rw-2rw^2}.$$
This proves the first equation in \eqref{eq}. We use this together with the equation for $a$ in
\eqref{eq:grad} and equation \eqref{beq} for $b$ to compute
\begin{eqnarray} \nonumber
\frac{az_0^2}{az_0^2+b}&=&\frac{az_0^8}{az_0^8+bz_0^6}\\ \nonumber
&=&\frac{\frac{rz_0^8}{2||z^2||^4}}{\frac{rz_0^8}{2||z^2||^4}+\frac{z_0^6}{2||z||^6}-
\frac{3rz_0^6}{4||z||^2||z^2||^2}}\\ \nonumber
&=&\frac{\frac{rw^2}{2v^2}}{\frac{rw^2}{2v^2}+\frac{1}{2v^3}-\frac{3rw}{4v^2}}\\ \nonumber
&=&\frac{rw^2}{rw^2+\frac{1}{v}-\frac{3rw}{2}}\\
&=&\frac{rw^2}{2}.\label{quot}
\end{eqnarray}
Dividing (\ref{e1}) by (\ref{e0}) and combining the result with (\ref{quot}) we get
$$||z||^2=z_0^2 \cdot \frac{I_1}{I_0}\bigg(-\frac{az_0^2}{az_0^2+b}\bigg)=
z_0^2 \cdot \frac{I_1}{I_0}\bigg(-\frac{rw^2}{2}\bigg).$$
Hence by definition of $v$ this can be rephrased as
$$v=\frac{I_1}{I_0}\bigg(-\frac{rw^2}{2}\bigg).$$
This proves the second equation in (\ref{eq}) and completes the proof in the case $r \neq 0$.
\\ \\
\textbf{Case\,2: } $r=0$.
\\ \\
In the case $r=0$ equation (\ref{eq}) becomes
\begin{equation}\label{simpl}
v=\frac{1}{2}=\frac{I_1(0)}{I_0(0)}.
\end{equation}
By equation~\eqref{be} in Appendix~\ref{sec:ell-int} we have
$$I_n(0)=\frac{(2n-1)!!\pi}{2^{n+1}n!},$$
where $(2n-1)!!$ equals $(2n-1)(2n-3)\cdots 1$ for $n\geq 1$ and $1$
for $n=0$. In particular,
$$I_0(0)=\frac{\pi}{2}, \qquad I_1(0)=\frac{\pi}{4}$$
and therefore the second equality in (\ref{simpl}) follows. It remains to check the first equation
in (\ref{simpl}) which says that
\begin{equation}\label{mit}
||z||^2=\frac{z(0)^2}{2}.
\end{equation}
However, for $r=0$ we have that $a=0$ and therefore $z$ is a solution of the ODE
$$z''=-bz,$$
which implies that up to scaling and time-reparametrization $z$ is given by the cosine function.
Now (\ref{mit}) follows from
$$\int_0^1 \cos^2(2\pi nt)dt=\frac{1}{2},\qquad n\in\N.$$
This finishes the proof of (\ref{eq}) in the case $r=0$ and the proposition follows. \hfill $\square$

%%%%%%%%%%%%%%%%%%%%%%%%%%%%%%%%%%%%%%%%%%%%%%%%%%%%%%%%%%%%%%%%
\subsection{Proof of nondegeneracy}\label{sec:nondeg}
%%%%%%%%%%%%%%%%%%%%%%%%%%%%%%%%%%%%%%%%%%%%%%%%%%%%%%%%%%%%%%%%

%marginpar{I left the old style with ``linearize the equation'' here EV.}
In this section we prove Theorem~\ref{nondeg}. 
Assume first that $r>0$.
Let $z$ be a critical point of $\mathcal{F}_r$
%of collision type,
i.e., a solution of the problem 
\begin{eqnarray}\label{eq:crit2}
z''&=&-bz-2az^3,\\ \nonumber
b&=&\frac{1}{2||z||^6}-\frac{3r}{4||z||^2\cdot ||z^2||^2},\\ \nonumber
a&=&\frac{r}{2||z^2||^4}.
\end{eqnarray}
%passing through the origin.
As in the previous section, after a time shift we may assume that $z$
attains its maximum $||z||_0=z(0)$ at $\tau=0$, hence $z'(0)=0$.
An element $\xi$ in the kernel of the Hessian of
$\mathcal{F}_r$ at $z$ is a solution of the linearized problem
\begin{equation}\label{prob}
\xi''=-b\xi-6az^2\xi-db(\xi)z-2da(\xi)z^3.
\end{equation}
In order to prove nondegeneracy we need to show that $\xi$ is a
constant multiple of $z'$. Note that 
$$da(\xi)=-\frac{4r\langle \xi,z^3\rangle}{||z^2||^6}, \qquad
db(\xi)=-\frac{3\langle \xi,z\rangle}{||z||^8}+\frac{3r\langle \xi,z\rangle}{2||z||^4\cdot ||z^2||^2}+\frac{3r\langle \xi,z^3\rangle}{||z||^2\cdot ||z^2||^4}.$$
Using (\ref{eq:crit2}) and (\ref{prob}) we obtain via integration by parts
\begin{eqnarray*}
-b\langle \xi,z\rangle-2a\langle \xi,z^3\rangle&=&\langle \xi,z''\rangle\\
&=&\langle \xi'',z \rangle\\
&=&-b\langle \xi,z\rangle-6a\langle \xi,z^3\rangle-db(\xi)||z||^2-2da(\xi)||z^2||^2
\end{eqnarray*}
and therefore
$$4a\langle \xi,z^3\rangle=-db(\xi)||z||^2-2da(\xi)||z^2||^2.$$
Plugging into this equation the formulas for $a$, $da$, and $db$, this becomes
$$\frac{2r\langle \xi,z^3\rangle}{||z^2||^4}=
\frac{3\langle \xi,z\rangle}{||z||^6}-\frac{3r\langle \xi,z\rangle}{2||z||^2\cdot ||z^2||^2}-\frac{3r\langle \xi,z^3\rangle}{||z^2||^4}+
\frac{8r\langle \xi,z^3\rangle}{||z^2||^4},$$
which simplifies to
\begin{equation}\label{non1}
\frac{r\langle \xi,z\rangle}{2||z||^2\cdot ||z^2||^2}-\frac{\langle \xi,z\rangle}{||z||^6}=\frac{r\langle \xi,z^3\rangle}{||z^2||^4}.
\end{equation}
In particular (recall our assumption $r>0$), we see from (\ref{non1}) that if $\langle \xi,z\rangle$ vanishes the same
has to hold for $\langle \xi, z^3\rangle$. 
\\ \\
We recall the quantities
$$v=\frac{||z||^2}{||z||_0^2}, \qquad w=\frac{||z||^2 \cdot ||z||_0^2}{||z^2||^2}$$
where $||z||_0=z(0)$. Their variations with respect to $\xi$ are given by
\begin{equation}\label{var}
\widehat{v}=\frac{2\langle \xi,z \rangle}{||z||_0^2}-\frac{2||z||^2 \xi_0}{||z||_0^3},
\qquad \widehat{w}=\frac{2\langle \xi,z\rangle ||z||_0^2+2||z||^2 ||z||_0\xi_0}{||z^2||^2}
-\frac{4\langle \xi, z^3\rangle ||z||^2||z||_0^2}{||z^2||^4}
\end{equation}
where $\xi_0:=\xi(0)$. From (\ref{eq}) we infer that
\begin{equation}\label{eqvar1}
\widehat{v}=\frac{d}{dw}\bigg(\frac{2}{4+3rw-2rw^2}\bigg)\widehat{w}
\end{equation}
and
\begin{equation}\label{eqvar2}
0=\frac{d}{dw}\Bigg(\big(4+3rw-2rw^2\big)\frac{I_1}{I_0}\bigg(-\frac{rw^2}{2}\bigg)\Bigg)
\widehat{w}.
\end{equation}
As explained in Appendix~\ref{sec:ell-int} in formula (\ref{ric2}), the quotient $\tfrac{I_1}{I_0}$ of
elliptic functions satisfies as a function of $m$ the Riccati differential equation
$$
\bigg(\frac{I_1}{I_0}\bigg)'=\frac{1}{2m(1-m)}-\frac{1}{m(1-m)}\frac{I_1}{I_0}+\frac{1}{2(1-m)}
\bigg(\frac{I_1}{I_0}\bigg)^2.
$$
Using this equation and (\ref{eq}) again we compute the derivative
in (\ref{eqvar2}) as follows:
\begin{eqnarray*}
& &\frac{d}{dw}\Bigg(\big(4+3rw-2rw^2\big)\frac{I_1}{I_0}\bigg(-\frac{rw^2}{2}\bigg)\Bigg)\\
&=&\big(3r-4rw\big)\frac{I_1}{I_0}\bigg(-\frac{rw^2}{2}\bigg)
-rw\big(4+3rw-2rw^2\big)\bigg(\frac{I_1}{I_0}\bigg)'\bigg(-\frac{rw^2}{2}\bigg)\\
&=&\frac{6r-8rw}{4+3rw-2rw^2}+\frac{2(4+3rw-2rw^2)}{w(2+rw^2)}-\frac{8}{w(2+rw^2)}\\
& &-\frac{4rw}{(2+rw^2)(4+3rw-2rw^2)}\\
&=&\frac{w(2+rw^2)(6r-8rw)+2(4+3rw-2rw^2)^2-8(4+3rw-2rw^2)-4rw^2}{w(2+rw^2)(4+3rw-2rw^2)}\\
&=&\frac{36rw-36rw^2-18r^2w^3+18r^2w^2}{w(2+rw^2)(4+3rw-2rw^2)}\\
&=&\frac{18r(1-w)(2+rw)}{(2+rw^2)(4+3rw-2rw^2)}.
\end{eqnarray*}
Since $r$ and $w$ are positive, we see from this formula that this derivative
vanishes only if $w=1$. In this case we obtain from (\ref{eq}) that
\begin{equation}\label{req}
\frac{2}{4+r}=\frac{I_1}{I_0}\bigg(-\frac{r}{2}\bigg).
\end{equation}
If we set $m=-\tfrac{r}{2}$ this amounts to the equation
$$1=(2-m)\frac{I_1}{I_0}(m).$$
By Lemma~\ref{mono} in Appendix~\ref{sec:ell-int} there are no solutions $m<0$ of this equation, and therefore
there are no solutions $r>0$ of equation (\ref{req}). Hence if $r>0$ we necessarily have
$\widehat{w}=0$, and therefore in view of (\ref{eqvar1}) as well $\widehat{v}=0$. From (\ref{var}) we
infer
\begin{equation}\label{non2}
\xi_0=\frac{\langle \xi,z \rangle ||z||_0}{||z||^2}, \qquad
\langle \xi, z^3\rangle=\frac{||z^2||^2}{2}\bigg(\frac{\langle \xi,z \rangle}{||z||^2}+\frac{\xi_0}{||z||_0}
\bigg),
\end{equation}
and therefore
$$\langle \xi, z^3\rangle=\frac{||z^2||^2}{||z||^2}\langle \xi,z \rangle.$$
Plugging this into (\ref{non1}) we obtain
$$\frac{r\langle \xi,z\rangle}{2||z||^2\cdot ||z^2||^2}-\frac{\langle \xi,z\rangle}{||z||^6}=\frac{r\langle \xi,z\rangle}{||z||^2\cdot ||z^2||^2},$$
or equivalently,
$$-\frac{\langle \xi,z\rangle}{||z||^6}=\frac{r\langle \xi,z\rangle}{2||z||^2\cdot ||z^2||^2}.$$
Since $r>0$, the two sides have opposite signs, therefore both sides have
to be zero and we obtain
$$\langle \xi,z\rangle=0.$$
Using (\ref{non1}) this implies
$$\langle \xi,z^3\rangle=0,$$
and therefore
$$da(\xi)=0, \qquad db(\xi)=0.$$
From $\langle \xi,z\rangle=0$ and the first equation in (\ref{non2}) we further conclude that 
$$\xi(0)=\xi_0=0.$$
Hence from (\ref{prob}) we see that $\xi$ is a solution of the ODE
\begin{equation}\label{eq:xi}
  \xi''=-b\xi-6az^2 \xi
\end{equation}
with $\xi(0)=0$. Applying the same reasoning to $z'$ in place of
$\xi$, we conclude that $z'$ also solves~\eqref{eq:xi} with $z'(0)=0$.
From~\eqref{eq:crit2} and~\eqref{eq:c-pos} we infer
$z''(0)=-b||z||_0-2a||z||_0^3<0$, so we can define
$$c:=\frac{\xi'(0)}{z''(0)}\in\R.$$
Then $\eta:=\xi-cz'$ solves~\eqref{eq:xi} with $\eta(0)=\eta'(0)=0$,
hence $\eta\equiv 0$ and $\xi=cz'$. 

This proves nondegeneracy for the case $r>0$. In the case
$r=0$ the functional $\mathcal{F}_0$ is just the functional for the regularized free fall
for which nondegeneracy can be checked directly, see \cite[Lemma\,3.6]{frauenfelder-weber}.
This finishes the proof of Theorem~\ref{nondeg}.

%%%%%%%%%%%%%%%%%%%%%%%%%%%%%%%%%%%%%%%%%%%%%%%%%%%%%%%%%%%%%%%%
\subsection{Uniqueness of symmetric critical points}\label{sec:uniqueness}
%%%%%%%%%%%%%%%%%%%%%%%%%%%%%%%%%%%%%%%%%%%%%%%%%%%%%%%%%%%%%%%%

%marginpar{NEW}
In this subsection we prove a uniqueness result for critical points of
$\FF_r$, $r\geq 0$.
To formulate the result, we introduce some terminology
from~\cite{cieliebak-frauenfelder-volkov}. A {\em symmetric critical point}
of $\FF_r$ is a smooth map $z:\R/2\Z\to\R$
satisfying the critical point equation~\eqref{eq:crit} and the
symmetry conditions
\begin{equation}\label{eq:symm}
   z(1+\tau)=-z(\tau)\text{ and } z(\tau) = z(1-\tau)\text{ for all }\tau.
\end{equation}
By the discussion at the beginning of this section and
Lemma~\ref{lem:LC}, $z$ has an odd number of zeroes in the
interval $[0,1)$ all of which are nondegenerate. 
Its Levi-Civit\'a transform $q:S^1\to\R$ has an odd number of zeroes
and satisfies $q(1-t)=q(t)$. Note that~\eqref{eq:symm} implies
$$
   z(0) = z'(1/2) = 0. 
$$
In particular, restriction to symmetric critical points removes
the translation invariance of the functional $\FF_r$. 

Symmetric critical points of $\FF_r$ may still be
nonunique because they may be multiply covered. Therefore, we restrict 
to symmetric critical points $z$ which are {\em simple}, i.e., of minimal period
$2$. This is equivalent to $z$ having zeroes precisely at integer
points $\tau\in\Z$, and critical points at $\tau\in1/2+\Z$. The
remaining ambiguity $z\mapsto -z$ can be removed by requiring $z$ to
be {\em normalized} by $z(\tau)>0$ for all $\tau\in(0,1)$. Using
Theorem~\ref{nondeg} we will prove 

\begin{cor}\label{cor:uniqueness}
For each $r\in[0,\infty)$, the functional $\FF_r$ has a unique
normalized simple symmetric critical point. This critical point is
nondegenerate of index zero. 
\end{cor}

%Proposition~\ref{prop:z-to-q} above and Lemma~D.2 of 
%\cite{cieliebak-frauenfelder-volkov} imply the following compactness statement.
%marginpar{This is no compactness statement.}
The proof is based on the following lemma. 

\begin{lemma}\label{lem:comp}
Let $z$ be a simple symmetric critical point of the functional
$\mathcal{F}_r$. Then its
$C^0$-norm $||z||_0$ satisfies
\begin{equation}\label{eq:bddbelow}
1\le ||z||_0
\end{equation}
and
\begin{equation}\label{eq:bddabove}
  ||z||_0\le \sqrt{2+(2r)^{1/3}}.
\end{equation}
\end{lemma}

\textbf{Proof: }
After a suitable time reparametrization, the Levi-Civita
transformation $q$ of $z$ satisfies the hypotheses of the loop $q_2$
in~\cite[Lemma~D.2]{cieliebak-frauenfelder-volkov}. 
This gives us the estimates
\begin{equation}\label{eq:est0}
  1\le ||q||_0
\end{equation}
%\begin{align*}
%&||q||_0\le 2\overline{q},\\
%&||q||_0\le 2+\frac{2r}{(2\overline{q})^2}  
%\end{align*}
%Therefore,
\begin{equation}\label{eq:est}
  ||q||_0\le \min\left\{2\overline{q},2+\frac{2r}{(2\overline{q})^2}\right\},
\end{equation}
where $\ol q>0$ is the average of $q$. 
Inequality~\eqref{eq:est0} via Levi-Civita transformation implies that 
$1\le ||z||_0$, proving~\eqref{eq:bddbelow}.

To get an estimate on the 
norm observe that for $\omega\ge 0$ and $x>0$ we have 
$$
\min\left\{x,\frac{\omega}{x^2}\right\}\le \omega^{1/3}.
$$
We use this for $\omega:=2r$ and $x:=2\overline{q}$ together with~\eqref{eq:est} to get
$$
  ||q||_0\le 2+(2r)^{1/3}.
$$
In view of $||q||_0=||z ||_0^2$ this proves~\eqref{eq:bddabove} completing
the proof of the lemma.
\hfill $\square$
\smallskip

\textbf{Proof of Corollary~\ref{cor:uniqueness}: }
Following Sections~6.2 and ~6.3 of~\cite{cieliebak-frauenfelder-volkov},
we introduce the Hilbert space of {\em symmetric loops} 
\begin{align*}
  H_{\rm sym}^2(S^1,\R) 
  := \bigl\{&z\in H^2(\R/2\Z,\R)\;\bigl| 
  -z(1+\tau)=z(\tau) = z(1-\tau)\text{ for all }\tau\bigr\}
\end{align*}
and its open subset
$$
  X := \{z\in H_{\rm sym}^2(S^1,\R) \mid z'(0)>0,\;z(\tau)>0\text{ for
    all }\tau\in(0,1)\}.
$$
We consider on $H_{\rm sym}^2(S^1,\R)$ the $L^2$-inner product $\la
z,w\ra = \int_0^1z(\tau)w(\tau)d\tau$ (which is an inner
product in view of the condition $z(1+\tau)=-z(\tau)$).
By the discussion at the beginning of this section, $\FF_r:X\to\R$ has
an $L^2$-gradient $\nabla\FF_r$ of class $C^1$.  
For $R>0$ consider the set
\begin{align*}
   \mathcal{Z}_R &:= \{(r,z)\mid \nabla\FF_r(z)=0\} 
   \subset [0,R]\times X
\end{align*}
with its induced topology. By the preceding discussion,
$\mathcal{Z}_R$ consists of pairs $(r,z)$ such that $z$ is a normalized simple
symmetric critical point of $\FF_r$.

Next we show that $\mathcal{Z}_R$ is compact. Recall that for
$(r,z)\in\mathcal{Z}_R$ the loop $z$ satisfies the ODE~\eqref{eq:crit}
with constants $a,b$ depending on $(r,z)$ given in~\eqref{eq:grad}. 
Thus the $C^0$-bound from Lemma~\ref{lem:comp}
and equation~\eqref{eq:crit} give a uniform $C^3$-bound on $z$ for
$(r,z)\in\mathcal{Z}_R$, and thus compactness of $\mathcal{Z}_R$ by
the Arzel\`a-Ascoli theorem, provided we have uniform bounds on the constants $a,b$. 
Note that by definition $a\ge 0$.
We will use the quantities $v,w>0$ defined in~\eqref{eq:vw}. 

\textbf{Case\,1:} $r=0$. In this case 
by definition $a=0$. Since $v=||z||^2/||z||_0^2$,
equation~\eqref{eq} yields $2||z||^2=||z||_0^2$.
Now~\eqref{beq} and~\eqref{eq:bddbelow}
imply boundedness of $b$.

\textbf{Case\,2:} $r>0$. In this case let us denote 
$\mathcal{A}:=rw^2>0$. Using this notation we can write
$$
\frac{1}{w}=\frac{\sqrt{r}}{\sqrt{\mathcal{A}}}\le 
\frac{\sqrt{R}}{\sqrt{\mathcal{A}}}.
$$
Equation~\eqref{eq} together with $v>0$ implies
$$
4+3\frac{1}{w}\mathcal{A}-2\mathcal{A} > 0.
$$
Together with the above upper bound on $\frac{1}{w}$
this gives us 
$$
4+3\sqrt{R}\sqrt{\mathcal{A}}-2\mathcal{A} > 0,
$$
which implies a uniform upper bound on $\mathcal{A}$. Therefore,
equation~\eqref{quot} implies 
\begin{equation}\label{eq:doubleineq}
0\le \frac{a||z||_0^2}{a||z||_0^2+b}\le \tilde C
\end{equation}
for some constant $\tilde C$ independent of $r$.
Hence the definition of $I_0$ yields
\begin{equation}\label{eq:bddI0}
  0<I_0(-\tilde C)\le I_0\left(-\frac{a||z||_0^2}{a||z||_0^2+b}\right)\le I_0(0) <\infty
\end{equation}
for all $r\in [0,R]$. 
Equation~\eqref{e0} (with $n=1$ since $z$ is a simple orbit)
then gives us constants $c,C>0$ with
\begin{equation}\label{eq:doubleineq2}
c\le a||z||_0^2+b\le C
\end{equation}
for all $r\in [0,R]$. The second inequality in~\eqref{eq:doubleineq2}  
together with with~\eqref{eq:doubleineq} 
gives us an upper bound on $a||z||_0^2$.
The latter together 
with~\eqref{eq:bddbelow} gives an
upper bound on $a$. The bounds are uniform with respect to 
$r\in [0,R]$. The double 
inequality~\eqref{eq:doubleineq2}
together with an upper bound on 
$a||z||_0^2$ gives us an upper bound on 
$|b|$. This concludes the proof of
compactness of $\mathcal{Z}_R$. 

On the other hand, according Theorem~\ref{nondeg} each critical point
of $\FF_r$ is nondegenerate. This implies that $\mathcal{Z}_R$ is
(as a preimage of a regular value) a compact $1$-dimensional submanifold transverse to the slices 
$\{r={\rm
  const}\}$. 
To see transversality, assume by contradiction
that $\mathcal{Z}_R$ is not transverse 
to the slice $\{r={\rm  const}\}$
at $(r^*,z^*)\in \mathcal{Z}_R$. 
Then there exists a local parametrization 
$$
(-\eps,\eps)\ni s\mapsto 
(r(s),z(s))\in \mathcal{Z}_R
$$
with $(r(0),z(0))=(r^*,z^*)$
and $\dds|_{s=0}r(s)=0$, whereas
$v:=\dds|_{s=0}z(s)\ne 0$.
We differentiate the critical point 
equation
$$
\nabla\FF_{r(s)}(z(s))=0
$$
at $s=0$ and use the chain rule to get 
$$
D\nabla\FF_{r^*}(z^*)v=0,
$$
contradicting nondegeneracy 
of $z^*$.
Now it is easy to see that for $r=0$ (describing the free fall) there
exists a unique normalized simple symmetric critical point $z_0$. By
the proof of Proposition D.4 in~\cite{cieliebak-frauenfelder-volkov},
the Hessian of $\FF_0$ at $z_0$ is positive definite, so $z_0$ is
nondegenerate of index zero. Therefore, $\mathcal{Z}_R$ intersects
each slice $\{r={\rm const}\}$ in a single point $z_r$, which is
nondegenerate of index zero as a critical point of $\FF_r$.
Since $R$ was arbitrary, this proves the corollary. 
\hfill $\square$

%%%%%%%%%%%%%%%%%%%%%%%%%%%%%%%%%%%%%%%%%%%%%%%%%%%%%%%%%%%%%%%%
\section{Frozen planet orbits for mean interaction}\label{sec:appl-frozen}
%%%%%%%%%%%%%%%%%%%%%%%%%%%%%%%%%%%%%%%%%%%%%%%%%%%%%%%%%%%%%%%%

%In this section we apply Theorem~\ref{nondeg} to show nondegeneracy of 
%critical points of the functional $\mathcal{B}_{av}$
%introduced in \cite{cieliebak-frauenfelder-volkov}.

%\subsection{Recalling $\mathcal{B}_{av}$.}
Let us recall the variational approach for frozen planet orbits in
helium from~\cite{cieliebak-frauenfelder-volkov}. For the mean
interaction between the two electrons, one considers the space
$$\mathcal{H}_{av}=\bigg\{(z_1,z_2) \in H^2(S^1,\mathbb{R}^2)\;\Bigl|\; ||z_1||>0,\,\,
||z_2||>0,\,\,\frac{||z_1^2||^2}{||z_1||^2}>\frac{||z_2^2||^2}{||z_2||^2}\bigg\}$$
and defines on it the functional
$\mathcal{B}_{av} \colon \mathcal{H}_{av} \to \mathbb{R}$
by
$$
\mathcal{B}_{av}(z_1,z_2):=2\sum_{i=1}^2\bigg(||z_i||^2\cdot ||z_i'||^2+\frac{1}{||z_i||^2}\bigg)-
\frac{||z_1||^2\cdot ||z_2||^2}{||z_1^2||^2\cdot ||z_2||^2-||z_2^2||^2\cdot ||z_1||^2}.
$$
Here $z_1$ and $z_2$ correspond to the Levi-Civit\`a regularizations
of the outer and inner electron, respectively.
We restrict to Sobolev class $H^2$ right away because on this space the 
$L^2$-gradient of $\mathcal{B}_{av}$ will be of class $C^1$. This
$L^2$-gradient is given by
\begin{eqnarray}\label{eq:critav}
\nabla\mathcal{B}_{av}(z_1,z_2)&=&\Bigl(-4||z_1||^2\mathcal{V}_1(z_1,z_2),
-4||z_2||^2\mathcal{V}_2
(z_1,z_2)\Bigr),
\\ \nonumber
\mathcal{V}_1(z_1,z_2)&=&-z_1''+a_1 z_1+b_1z_1^3,\qquad
\mathcal{V}_2(z_1,z_2)=-z_2''+a_2 z_2+b_2z_2^3,
\\ \nonumber
a_1&=&\frac{||z_1'||^2}{||z_1||^2}-\frac{1}{||z_1||^6}-\frac{||z_2||^4 \cdot ||z_1^2||^2}{2||z_1||^2 \cdot\big(||z_1^2||^2\cdot ||z_2||^2
-||z_2^2||^2 \cdot ||z_1||^2\big)^2},
\\ \nonumber
b_1&=&\frac{||z_2||^4}{\big(||z_1^2||^2\cdot ||z_2||^2-||z_2^2||^2 \cdot ||z_1||^2\big)^2},\\
\nonumber
a_2&=&\frac{||z_2'||^2}{||z_2||^2}-\frac{1}{||z_2||^6}+\frac{||z_1||^4 \cdot ||z_2^2||^2}{2||z_2||^2 \cdot\big(||z_1^2||^2\cdot ||z_2||^2
-||z_2^2||^2 \cdot ||z_1||^2\big)^2},
\\ \nonumber
b_2&=&-\frac{||z_1||^4}{\big(||z_1^2||^2\cdot ||z_2||^2-||z_2^2||^2 \cdot ||z_1||^2\big)^2}.
\end{eqnarray}
It was shown in~\cite{cieliebak-frauenfelder-volkov} that
$\nabla\mathcal{B}_{av}$ is of class $C^1$; we refer to 
its derivative as the {\em Hessian} of $\mathcal{B}_{av}$. 
According to Lemma~3.1 in \cite{frauenfelder1} or Lemma~D2
of~\cite{cieliebak-frauenfelder-volkov}, for each critical point
$(z_1,z_2)$ of $\BB_{av}$ the first component $z_1$ is constant. 
Thus the $S^1$-action by time shift acts trivially on the first
component and 
(by a similar argument as for 
$\FF_r$) nontrivially on the second one, leading to a
$1$-dimensional subspace of the kernel of the Hessian. 
By analogy with $\mathcal{F}_r$, we say that 
a critical point of $\mathcal{B}_{av}$ is 
{\em nondegenerate} if the nullity of its Hessian is precisely $1$.
Now we can formulate the main result of this section, whose proof will
occupy the rest of this section. 

\begin{thm}\label{thm:nondeg-av}
All critical points of $\mathcal{B}_{av}$ are nondegenerate.
\end{thm}

%%%
\subsection{Relating frozen planet orbits to critical points of $\mathcal{F}_\rho$}
%%%

We introduce the numerical parameters
$$
\rho:=(\sqrt{2}-1)^2,\qquad \alpha:=\frac{\sqrt{2}-1}{\sqrt{2}}
$$
and observe that 
\begin{equation*}%\label{eq:aless1}
0<\alpha<1.
\end{equation*}
We introduce the following constant for any 
$z\in H^2_*(S^1,\R)$:
$$
c(z):=\alpha^{-1/2}\frac{||z^2||}{||z||}.
$$

\begin{lemma}
The frozen functional $\mathcal{F}_\rho$ is related to the functional
$\mathcal{B}_{av}$ by
\begin{equation}\label{twofu}
  \mathcal{F}_\rho(z)=\mathcal{B}_{av}\bigl(c(z),z\bigr).
\end{equation}
\end{lemma}

{\bf Proof: }
This follows from the following computation:
\begin{eqnarray*}
& &\mathcal{B}_{av}\Bigg(\frac{2^{1/4}\cdot||z^2||}{\sqrt{\sqrt{2}-1}\cdot||z||},z\Bigg)\\
&=&2||z||^2||z'||^2+\frac{2}{||z||^2}+\frac{2(\sqrt{2}-1)\cdot||z||^2}{\sqrt{2}\cdot ||z^2||^2}-\frac{||z||^2}{\frac{\sqrt{2}||z^2||^2}{(\sqrt{2}-1)||z||^2}||z||^2-||z^2||^2}\\
&=&2||z||^2||z'||^2+\frac{2}{||z||^2}+\frac{\sqrt{2}(\sqrt{2}-1)\cdot||z||^2}{||z^2||^2}-\frac{||z||^2}{\Big(\frac{\sqrt{2}}{\sqrt{2}-1}-1\Big)||z^2||^2}\\
&=&2||z||^2||z'||^2+\frac{2}{||z||^2}+\frac{\sqrt{2}(\sqrt{2}-1)\cdot||z||^2}{||z^2||^2}-\frac{(\sqrt{2}-1)||z||^2}{||z^2||^2}\\
&=&2||z||^2||z'||^2+\frac{2}{||z||^2}+\big(\sqrt{2}-1\big)^2\frac{||z||^2}{||z^2||^2}\\
&=&\mathcal{F}_\rho(z).
\end{eqnarray*}
\hfill$\square$

Let us illustrate the significance of the constant $c(z)$
from another angle. Let
$$
H^2_{const}:=\{(z_1,z_2)\in H^2(S^1,\R^2)\mid z_1=const\}
$$
be the subspace of functions with constant first 
component and set
\begin{align}\label{eq:avconst}
  \mathcal{H}_{av}^{const}:=\mathcal{H}_{av}\cap H^2_{const}=&\bigg\{(z_1,z_2)\in H^2(S^1,\R^2)\\ \nonumber
&\ \ \ \Bigl|\; z_1=const,\,\, ||z_2||>0,\,\, z_1^2>
\frac{||z_2^2||^2}{||z_2||^2}\bigg\}.
\end{align}
Consider the following codimension $1$ Hilbert submanifold
of $H^2_{const}$:
$$
graph(c):=\{(c(z),z)\in H^2_{const}\mid z\in H^2_*(S^1,\R)\}.
$$
Then the definition of $c(z)$ and $\alpha<1$ imply
$$
graph(c)\subset\mathcal{H}_{av}^{const}.
$$

\begin{lemma}\label{lem:vanish1}
Assume that 
$(z_1,z_2)\in \mathcal{H}_{av}^{const}$. Then the equation
$$
\mathcal{V}_1(z_1,z_2)=0
$$
is equivalent to the equation
$$
z_1=\pm c(z_2).
$$
Moreover, if this is the case, then
$$
a_1=-\frac{2}{z_1^6},\,\,\,\, b_1=\frac{2}{z_1^8}.
$$
\end{lemma}

\textbf{Proof: } For constant $z_1$, the functions 
$a_1$ and $b_1$ that enter the formula for 
$\mathcal{V}_1$ simplify to
\begin{equation}\label{fr0}
a_1=-\frac{1}{z_1^6}-\frac{||z_2||^4z_1^2}{2\big(||z_2||^2z_1^4-||z_2^2||^2z_1^2\big)^2},
\qquad b_1=\frac{||z_2||^4}{\big(||z_2||^2z_1^4-||z_2^2||^2z_1^2\big)^2}.
\end{equation}
Therefore, $\mathcal{V}_1$ rewrites as follows:
$$
 \mathcal{V}_1(z_1,z_2)=a_1z_1+b_1z_1^3=-\frac{1}{z_1^5}+\frac{||z_2||^4z_1^3}{2\big(||z_2||^2z_1^4-||z_2^2||^2z_1^2\big)^2}.
$$
Thus, $\mathcal{V}_1(z_1,z_2)=0$ is equivalent to
\begin{equation}\label{fr1}
2\big(||z_2||^2z_1^4-||z_2^2||^2z_1^2\big)^2=||z_2||^4z_1^8.
\end{equation}
The defining inequality in~\eqref{eq:avconst} implies that~(\ref{fr1}) 
is equivalent to
$$
\sqrt{2}\big(||z_2||^2z_1^4-||z_2^2||^2z_1^2\big)=
||z_2||^2z_1^4.
$$
Since $z_1\ne 0$ by~\eqref{eq:avconst},
the latter is equivalent to
$$
||z_2^2||^2=\frac{\sqrt{2}-1}{\sqrt{2}}||z_2||^2z_1^2,
$$
in other words to
$$
z_1=\pm\frac{2^{1/4}\cdot||z_2^2||}{\sqrt{\sqrt{2}-1}\cdot||z_2||}.
$$
This shows the first statement. 
The second statement follows by substituting~\eqref{fr1}
in~\eqref{fr0}.
\hfill $\square$

Taking the derivative of equation~\eqref{twofu} with respect to $z$ yields
\begin{equation*}
D\mathcal{F}_\rho(z)=D_1\mathcal{B}_{av}(c(z),z) Dc(z) + D_2\mathcal{B}_{av}(c(z),z)
= D_2\mathcal{B}_{av}(c(z),z).
\end{equation*}
Here $D_1$ and $D_2$ denotes the derivative with respect to the first
and second component, respectively, and
$D_1\mathcal{B}_{av}(c(z),z)=0$ by Lemma~\ref{lem:vanish1}. 
As a consequence, we get the following relation 
between the $L^2$-gradients of   
$\mathcal{F}_\rho(z)$ and $\mathcal{B}_{av}$:
\begin{equation}\label{eq:relgrad2}
\nabla \mathcal{F}_\rho(z)= \mathcal{V}_2(c(z),z).
\end{equation}
Here we have dropped the factor $-4||z||^2$ from
$\nabla\mathcal{F}_\rho(z)$ which has no relevance for the subsequent discussion.
Taking another derivative with respect to $z$, we get
\begin{equation}\label{eq:UrsDiff}
D\nabla \mathcal{F}_\rho(z)\xi=
D_1\mathcal{V}_2(c(z),z)Dc(z)\xi+
D_2\mathcal{V}_2(c(z),z)\xi
\end{equation}
for any $\xi\in H^2(S^1,\R)$. 

\subsection{Proof of Theorem~\ref{thm:nondeg-av} modulo two key lemmas}\label{ss:modkeylemmas}

Let $(z_1,z_2)$ be a critical point of $\mathcal{B}_{av}$. 
Recall from~\cite{frauenfelder1,cieliebak-frauenfelder-volkov} 
that $z_1$ is constant. Therefore, Lemma~\ref{lem:vanish1} implies that 
$z_1=\pm c(z_2)$. In the following we assume $z_1>0$ (the case $z_1<0$
being analogous), so that
$$
(z_1,z_2)\in graph(c). 
$$
Equation~\eqref{eq:relgrad2} with $z=z_2$ implies that
$z_2$ is a critical point of $\mathcal{F}_\rho$.

Let $\xi=(\xi_1,\xi_2)\in H^2(S^1,\R^2)$ lie in the 
kernel of the Hessian at $(z_1,z_2)$, that is 
\begin{equation}\label{eq:hess1}
D_1\mathcal{V}_1(z_1,z_2)\xi_1+
D_2\mathcal{V}_1(z_1,z_2)\xi_2=0
\end{equation}
and 
\begin{equation}\label{eq:hess2}
D_1\mathcal{V}_2(z_1,z_2)\xi_1+
D_2\mathcal{V}_2(z_1,z_2)\xi_2=0.
\end{equation}

We need the following two lemmas.

\begin{lemma}\label{lem:const}
Assume that 
$(z_1,z_2)\in graph(c)$ and
$(\xi_1,\xi_2)\in H^2(S^1,\R^2)$ are such that 
equation~\eqref{eq:hess1} is satisfied. Then
$\xi_1$ is constant, that is $(\xi_1,\xi_2)\in H^2_{const}$.
\end{lemma}

\begin{lemma}\label{lem:xi1xi2}
Assume that $(z_1,z_2)\in graph(c)$
and $(\xi_1,\xi_2)\in H^2_{const}$
are such that equation~\eqref{eq:hess1} is satisfied. Then
$$
\xi_1=Dc(z_2)\xi_2,
$$
that is $(\xi_1,\xi_2)\in T_{(z_1,z_2)}graph(c)$.
\end{lemma}

Assuming these two lemmas, set $z=z_2$,
$c(z)=z_1$ and $\xi=\xi_2$ 
in equation~\eqref{eq:UrsDiff} and
solve it for 
$D_2\mathcal{V}_2(z_1,z_2)\xi_2$ to get 
$$
D_2\mathcal{V}_2(z_1,z_2)\xi_2=D_2\nabla
\mathcal{F}_\rho(z_2)\xi_2-D_1
\mathcal{V}_2(z_1,z_2)
Dc(z_2)\xi_2.
$$
The latter allows us to rewrite 
equation~\eqref{eq:hess2} as 
\begin{equation}\label{eq:UrsKey}
D\nabla\mathcal{F}_\rho(z_2)\xi_2=D_1\mathcal{V}_2(z_1,z_2)(Dc(z_2)\xi_2-\xi_1).
\end{equation}
Apply Lemma~\ref{lem:xi1xi2} to get
$$
D\nabla\mathcal{F}_\rho(z_2)\xi_2=0.
$$
Since $z_2$ is a critical point of 
$\mathcal{F}_\rho$, it is nondegenerate by Theorem~\ref{nondeg}, so
$\xi_2$ belongs to a subspace of dimension $1$. This 
together with Lemma~\ref{lem:xi1xi2} shows that 
$(\xi_1,\xi_2)$ belongs to a 
subspace of dimension $1$, giving us nondegeneracy of
$(z_1,z_2)$ as a critical point of $\mathcal{B}_{av}$.
This proves Theorem~\ref{thm:nondeg-av} modulo the two lemmas above.

\subsection{Proof of Lemma~\ref{lem:const}}

Observe that equation~\eqref{eq:hess1} is 
equivalent to 
$$
\xi_1''=a_1\xi_1+3b_1z_1^2\xi_1+\Bigl(da_1(\xi_1,\xi_2)z_1+db_1(\xi_1,\xi_2)
z_1^3\Bigr).
$$
Note that the term in the round brackets is a function constant in time, which we abbreviate
as $const$. So, we can rewrite it as
$$
\xi_1''=(a_1+3b_1z_1^2)\xi_1+const.
$$
Inserting for $a_1,b_1$ the simplified expressions from
Lemma~\ref{lem:vanish1}, we see that
$$
a_1+3b_1z_1^2=-\frac{2}{z_1^6}+\frac{6z_1^2}{z_1^8}=
\frac{4}{z_1^6}>0
$$
is a positive constant. Now we argue by contradiction.

Suppose that $\xi_1$ is not constant and let $(\tau_1,\tau_2)$ be a maximal
interval on which $\xi_2'>0$. 
%Let $t_{min}$ and $t_{max}$ be the times when 
%$\xi_1$ reaches minimum respectively maximum.
%Then acceleration at $t_{min}$ must be nonnegative, at $t_{max}$ nonpositive. Note that
Then $\xi_1'(\tau_1)=\xi_1'(\tau_2)=0$, $\xi_1''(\tau_1)\geq 0$, and
$\xi_1''(\tau_2)\leq 0$.  
Since the coefficient in front of $\xi_1$ is constant in time and
positive, we get the following chain of inequalities for $\tau_1\leq
\tau\leq \tau_2$:
$$
0\le \xi_1''(\tau_1)\le\xi_1''(\tau)\le\xi_1''(\tau_2)\le 0.
$$
So we must have equalities everywhere, $\xi''(\tau)=0$, and
therefore $\xi'(\tau)$ for all $\tau\in(\tau_1,\tau_2)$, contradicting
our hypothesis.

\subsection{Proof of Lemma~\ref{lem:xi1xi2}.}

{\bf Step 1. }
Let $z=(z_1,z_2)\in graph(c)$
and $\xi=(\xi_1,\xi_2)\in H^2_{const}$ satisfy equation~\eqref{eq:hess1}.
It follows from equation~\eqref{eq:critav} that the restriction 
$$
\mathcal{W}:=\mathcal{V}_1|_{H^2_{const}}
$$
lands in the space of constant functions and can, therefore, be considered as a map
$$
H^2_{const}\longrightarrow \R.
$$
Since both $z_1$ and $\xi_1$ are constant in $t$,
we can assume that we vary $z_1$ within the space 
of constant functions when computing the derivative $D_1\mathcal{V}_1$
with respect to $z_1$.  
Therefore, equation~\eqref{eq:hess1} is equivalent to 
\begin{equation}\label{eq:hess1const}
D\mathcal{W}(z)\xi=D_1\mathcal{W}(z_1,z_2)\xi_1+
D_2\mathcal{W}(z_1,z_2)\xi_2=0.
\end{equation}
In this terminology, Lemma~\ref{lem:vanish1} is equivalent to
$$
\mathcal{H}_{av}\cap\{\mathcal{W}=0\}=graph(c).
$$
Assume for the moment that
\begin{equation}\label{eq:reg}
D\mathcal{W}(z)\ne 0.
\end{equation}
Then
$$ 
 T_zgraph(c)=\ker D\mathcal{W}(z),
$$
since the inclusion 
$T_zgraph(c)\subset\ker D\mathcal{W}(z)$
is obvious and both spaces have codimension $1$
in $H^2_{const}$. This shows the lemma 
modulo~\eqref{eq:reg}, which
is an immediate consequence of the next step.

{\bf Step 2. }
For any $z=(z_1,z_2)\in graph(c)$, we have 
\begin{equation}\label{eq:nonvan}
D_1\mathcal{W}(z_1,z_2)\ne 0.
\end{equation}
This is achieved by a brute force computation. 
Since $z_1$ is constant in $t$, we set $z_1''=0$
in formula~\eqref{eq:critav} to get
\begin{equation}\label{eq:nu11}
\mathcal{W}(z_1,z_2)=a_1z_1+b_1z_1^3,
\end{equation}
where $a_1,b_1$ are given by~\eqref{fr0} to be
$$
a_1=-\frac{1}{z_1^6}-\frac{z_1^2}{2}b_1,\quad b_1=\frac{||z_2||^4}{P(z_1)^2}
$$
with the polynomial
$$
P(z_1):=||z_2||^2z_1^4-||z_2^2||^2z_1^2.
$$
We rewrite~\eqref{eq:nu11} as
\begin{equation}\label{eq:nu12}
\mathcal{W}(z_1,z_2)=-\left(\frac{1}{z_1^6}+\frac{z_1^2}{2}b_1\right)z_1+b_1z_1^3=-\frac{1}{z_1^5}
+\frac{z_1^3}{2}b_1
\end{equation} and
compute the derivative with respect to $z_1$: 
$$
D_1\mathcal{W}(z_1,z_2)=\frac{5}{z_1^6}+\frac{3}{2}z_1^2b_1+\frac{1}{2}z_1^3D_1b_1(z_1).
$$
Now we differentiate $b_1$:
$$
D_1b_1(z_1)=-2||z_2||^4\frac{P'(z_1)}{P(z_2)^3}=-2\frac{P'(z_1)}{P(z_1)}b_1
$$
Therefore,
\begin{equation}\label{eq:Dnu}
D_1\mathcal{W}(z_1,z_2)=\frac{5}{z_1^6}+\left(\frac{3}{2}z_1^2-z_1^3\frac{P'(z_1)}{P(z_1)}\right)b_1.
\end{equation}
Now we manipulate $P(z_1)$ and $P'(z_1)$
\begin{align*}
        P(z_1) =&||z_2||^2z_1^4-||z_2^2||^2z_1^2\\
  \stackrel{(*)}{=}&||z_2||^2z_1^4-\alpha ||z_2||^2z_1^4\\
       =&(1-\alpha )||z_2||^2z_1^4.
\end{align*}
\begin{align*}
P'(z_1)=&4||z_2||^2z_1^3-2||z_2^2||^2z_1\\
       =&z_1\bigl(4||z_2||^2z_1^2-2||z_2^2||^2\bigr)\\
       \stackrel{(*)}{=}&
      z_1\bigl(4||z_2||^2z_1^2-2\alpha ||z_2||^2z_1^2\bigr)\\
      =&2(2-\alpha)||z_2||^2z_1^3.
\end{align*}
Here the equalities marked with $(*)$ use $z_1=c(z_2)$
in the form
$$
||z_2^2||^2=\alpha ||z_2||^2z_1^2.
$$
Now we set
$$
X:=P(z_1)^3z_1^6D_1\mathcal{W}(z_1,z_2)
$$
and continue from \eqref{eq:Dnu}:
\begin{align*}
X=&5P(z_1)^3+\frac{3}{2}||z_2||^4z_1^8P(z_1)-||z_2||^4z_1^9P'(z_1)\\
=&\Bigl(5(1-\alpha)^3+\frac{3}{2}(1-\alpha)
-2(2-\alpha)\Bigr)||z_2||^6z_1^{12}.
\end{align*}
We compute the numerical coefficient 
$$
K:=5(1-\alpha)^3+\frac{3}{2}(1-\alpha)-2(2-\alpha)
$$
in front of
$||z_2||^6z_1^{12}\xi_1$.
With
$$
k:=1-\alpha=\frac{1}{\sqrt{2}}
$$
we compute:
\begin{align*}
K+2\alpha
=&5(1-\alpha)^3+\frac{3}{2}
(1-\alpha)-2(2-2\alpha)\\
=&5k^3+\frac{3}{2}k-4k\\
%    =&5k^3-\frac{5}{2}k\\
    =&5k(k^2-\frac{1}{2})
    =0.
\end{align*}
This shows $K=-2\alpha\ne 0$, which implies $X\neq 0$ and therefore
proves equation~\eqref{eq:nonvan}.

%%%%%%%%%%%%%%%%%%%%%%%%%%%%%%%%%%%%%%%%%%%%%%%%%%%%%%%%%%%%%%%%
\subsection{Uniqueness of symmetric frozen planet orbits}\label{sec:uniqueness-mean}
%%%%%%%%%%%%%%%%%%%%%%%%%%%%%%%%%%%%%%%%%%%%%%%%%%%%%%%%%%%%%%%%

The discussion in the preceding subsections shows that critical points
of $\BB_0=\BB_{av}$ (i.e., frozen planet orbits for the mean
interaction) are in one-to-one correspondence with critical
points of $\FF_\rho$ for $\rho=(\sqrt{2}-1)^2$. Moreover, this
correspondence preserves the indices and nullities of the critical points. 
By Corollary~\ref{cor:uniqueness}, the functional $\FF_\rho$ has a unique
normalized simple symmetric critical point, which is nondegenerate of
index zero. Hence the same holds for the corresponding normalized
simple symmetric critical point of $\BB_{av}$ and we conclude the
following result, which corresponds to Theorem A in the Introduction.

\begin{cor}\label{cor:uniqueness-mean}
The unique normalized simple symmetric frozen planet orbit for
the mean intersection functional $\BB_{av}$ is nondegenerate of
Morse index $0$.
\hfill$\square$
\end{cor}

%%%%%%%%%%%%%%%%%%%%%%%%%%%%%%%%%%%%%%%%%%%%%%%%%%%%%%%%%%%%%%%%%%%%%%%%%%
\section{Determinant lines of self-adjoint Fredholm operators}\label{sec:det}
%%%%%%%%%%%%%%%%%%%%%%%%%%%%%%%%%%%%%%%%%%%%%%%%%%%%%%%%%%%%%%%%%%%%%%%%%%

%%%
%\subsection{Index zero symmetric Fredholm operators.}\label{ss:ss}
%%%

The goal of this section is to show that under certain conditions the
restriction of the determinant bundle to the space of symmetric index
zero Fredholm operators bounded from below (or above) is trivial. We
begin by describing the general Fredholm setting. Then we will
describe the more specific Hilbert space setting and state the main
result. 
\smallskip

%%%
\subsection{The determinant line bundle}
%%%

For two Banach spaces $X$ and $Y$ we denote the space of continuous
linear maps from $X$ to $Y$ by $\LL(X,Y)$, the subspace of Fredholm
operators from $X$ to $Y$ by $\FF(X,Y)$, and the set of surjective Fredholm 
operators by $\FF^*(X,Y)$. Recall from \cite{salamon} that the determinants
$$
   \det(D) = \Lambda^{\rm max}(\ker D^*)\otimes\Lambda^{\rm max}(\ker D)
$$
for any $D\in \FF(X,Y)$ give rise to a real line bundle, the {\em determinant line bundle}
$$
   \det \to \FF(X,Y). 
$$
We describe the bundle structure for $\det$ following
\cite{salamon}. For this recall the bundle of kernels  
over $\FF^*(X,Y)$ and note that the restriction 
$\det|_{\FF^*(X,Y)}$ is just the top exterior power
of the kernel bundle over $\FF^*(X,Y)$. 
Let now $T\in \FF(X,Y)$ be a Fredholm operator and
set $N:=\dim\coker T$. Let 
$$
  \Phi:\R^N\longrightarrow Y
$$
be an isomorphism onto a direct complement to $\im T$.
Then the stabilized operator 
$$
  D\oplus \Phi:X\oplus \R^N\longrightarrow Y,\quad (x,\zeta)\mapsto D(x)+\Phi(\zeta)
$$
is surjective for $D:=T$ by construction. Therefore, 
it is surjective for any $D$ in a small 
open neighbourhood $U_{T,\Phi}$ of $T$. 
The idea is to construct a certain fiberwise linear bijection 
$\iota_\Phi$ between the restriction of the determinant bundle to
$U_{T,\Phi}$ and its restriction to the image of $U_{T,\Phi}$ under
stabilization. We then declare that  
$\iota_\Phi$ be a bundle isomorphism. To define $\iota_\Phi$
pick any $D\in U_{T,\Phi}$ and set $k:=\dim\ker D$, 
$l:=\dim\coker D$. Then 
$$
  \dim\ker(D\oplus\Phi)=\ind(D\oplus\Phi)=(k-l)+N
$$ 
and 
$$
  \ker D\times \{0\}\subset \ker(D\oplus\Phi)
$$
is a subspace of dimension $k$. In particular,
\begin{equation}\label{eq:dimineq}
  k=\dim\ker D\le\dim\ker(D\oplus\Phi)=N+(k-l).
\end{equation}
A complement to $\ker D\times \{0\}$ in $\ker(D\oplus\Phi)$ can be
described as follows. Pick $N-l$ linearly independent vectors 
$\{\zeta_j\}_{j=l+1}^N$ in $\R^N$ and $N-l$ vectors
$\{\xi_j\}_{j=l+1}^N$ in $X$ subject to 
$$
  D\xi_j+\Phi\zeta_j=0,\quad j=l+1,\dots,N.
$$
Then the collection $\{(\xi_j,\zeta_j)\}_{j=l+1}^N$ spans the 
desired complement. We complete $\{\zeta_j\}_{j=l+1}^N$ 
to a basis $\{\zeta_j\}_{j=1}^N$ of $\R^N$. Now for 
any 
$$
\theta=(y_l^*\wedge\dots\wedge y_1^*)\otimes
(x_1\wedge\dots \wedge x_k)\in \det(D) 
$$
define
$$
\iota_\Phi(D,\theta):=
(-1)^{kl}\frac{\det\la y_j^*,\Phi\zeta_{i}\ra_{i,j=1,\dots,l}}
{\det(\zeta_1,\dots,\zeta_N)}
(x_1,0)\wedge\dots\wedge(x_k,0)\wedge 
(\xi_{l+1},\zeta_{l+1})\wedge\dots\wedge 
(\xi_N,\zeta_N).
$$
It is shown in \cite{salamon} that the map 
$\iota_\Phi$ is independent of the choices made and the collection 
of all maps $\iota_\Phi$ does indeed define a bundle 
structure on the determinant bundle over $\FF(X,Y)$.

%%%
\subsection{Self-adjoint Fredholm operators}
%%%

Let now $(F,\la\ ,\ \ra_F)$ be a real Hilbert space and $E\subset F$ a
dense linear subspace. Recall that a linear map 
$$
T:E\longrightarrow F
$$
is called {\em symmetric} if  
$$
   \la Tx,y\ra_F = \la x,Ty\ra_F \quad \text{for all }x,y\in E,
$$
and {\em self-adjoint} if in addition for each $y\in F$ the existence
of a constant $C_y$ such that  
\begin{equation}\label{eq:selfadj}
\la Tx,y\ra_F\leq C_y\|x\|_F
\end{equation} 
for all $x\in E$ implies that $y\in E$. 

Assume now that $(E,\la\ ,\ \ra_E)$ is itself a Hilbert space, so we
can talk about the space of Fredholm operators $\FF(E,F)$. 

\begin{lemma}\label{lem:selfadj}
In the setting above, a Fredholm operator $T\in\FF(E,F)$ is
self-adjoint if and only if it is symmetric and of index zero. 
Moreover, in this case we have the equalities
\begin{equation}\label{eq:kerim}
  \ker T=(\im T)^{\perp_F},\qquad (\ker T)^{\perp_F}=\im T.
\end{equation}
\end{lemma}

\textbf{Proof: }
Suppose first that $T$ is self-adjoint. Then it is in particular
symmetric and we obtain the inclusion $\ker T\subset(\im T)^{\perp_F}$.
On the other hand, consider $y\in(\im T)^{\perp_F}$. Then $\la
Tx,y\ra =0$ for all $x\in E$, so $y$ satisfies~\eqref{eq:selfadj} with
$C_y=0$. Since $T$ is self-adjoint, this implies $y\in E$, and
symmetry of $T$ gives $\la x,Ty\ra=0$ for all $x\in E$. By density of
$E\subset F$ this implies $Ty=0$, so we have shown $\ker T=(\im
T)^{\perp_F}$. This proves $\ind\,T=0$ and the first equality, and the
second equality follows from the first one by taking orthogonal
complements. 

Suppose now that $T$ is symmetric and of index zero. 
Then the inclusion $\ker T\subset (\im T)^{\perp F}$ has to be an
equality because both spaces have the same finite dimension, so the
equalities~\eqref{eq:kerim} hold. 
To prove self-adjointness, assume first that $T$ is surjective. 
Let $y\in F$ satisfy~\eqref{eq:selfadj}. By the Riesz 
representation theorem in $F$ for the functional 
$x\mapsto \la Tx,y\ra_F$ and surjectivity of
$T$, there exists $z\in E$ with 
$$
\la Tx,y\ra_F=\la x,Tz\ra_F
$$
for all $x\in E$. We apply symmetry of $T$
to the right hand side of the last displayed equation to get 
$$
  \la Tx,y-z\ra_F=0
$$
for all $x\in E$. Since $\im T=F$, this implies $y-z=0$ and thus $y=z\in E$.
This proves self-adjointness in the case that $T$ is surjective.
%Let now $T\in \FF_s(E,F)$ be arbitrary.
If $T$ is not surjective, then we replace the triple $(F,E,T)$ by the triple
$$
  (F_1:=\im T,\,E_1:=\im T\cap E,\,T_1:=T|_{\im T\cap E}).
$$
Since the codimension of $\im T$ is finite, the linear subspace $E_1$
is dense in $F_1$ (because $\im T$ possesses a complement which is
contained in $E$).  
Equations~\eqref{eq:kerim} imply that $T_1$ is surjective, so by the
discussion above $T_1$ is self-adjoint. Therefore, the original
operator $T$ is self-adjoint.
\hfill $\square$

We denote the space of operators as in Lemma~\ref{lem:selfadj} by
$$
  \FF_s(E,F):=\{T\in \FF(E,F)\mid \ind\,T=0\,\,
\text{and}\,\, T\,\, \text{is symmetric}\}.
$$
Let $T\in\FF_s(E,f)$. By Lemma~\ref{lem:selfadj} and its proof, $\im
T\cap E$ is a complement to $\ker T$ in $E$. So the restriction of $T$
to $\im T\cap E$ defines an isomorphism
$$
  \wt T:= T|_{\im T\cap E}:(\im T\cap
  E,\la\ ,\ \ra_E)\stackrel{\cong}\longrightarrow (\im T,\la\ ,\ \ra_F).  
$$
Assume now in addition that the inclusion $E\into F$ is compact. Then
the inverse of $\wt T$ can be viewed as a compact operator
$$
  \wt T^{-1}:(\im T,\la\ ,\ \ra_F) \to (\im T\cap
  E,\la\ ,\ \ra_E)\into (\im T,\la\ ,\ \ra_F).
$$
As such, the spectrum of $\wt T^{-1}$ consists only of eigenvalues
which are uniformly bounded and can accumulate only at zero, and all
nonzero eigenvalues have finite multiplicity (see e.g.~\cite{conway}).
Since the eigenvalue zero of $T$ has only finite multiplicity by the
Fredholm hypothesis, we obtain 

\begin{lemma}\label{lem:spectrum}
Consider the setting above and assume that the inclusion $E\into F$ is
compact. Then for $T\in\FF_s(E,F)$ the spectrum $\sigma(T)\subset\R$
is discrete and consists only of eigenvalues of finite multiplicity.
\hfill$\square$
\end{lemma}

For $T$ as in the preceding lemma and a measurable subset $Z\subset\R$ let 
$$
  E_Z^T:F\to F
$$
denote the projection onto the sum of eigenspaces to eigenvalues in
$Z$ (this is a very special case of a {\em spectral
  measure}~\cite{conway}), so we get the spectral decomposition 
$$
  T = \bigoplus_{\lambda\in\sigma(A)}\lambda E_{\lambda}^T. 
$$
We conclude this subsection with an important lemma about the
continuity of the spectrum and of the spectral measure.  

\begin{lemma}
Consider the setting as in Lemma~\ref{lem:spectrum} and $T \in
\FF_s(E,F)$.

(a) Let $I\subset \R$ be a compact interval such that $T$ has no
eigenvalues in $I$. Then there exists an open neighbourhood $W_T$ of
$T$ in $\FF_s(E,F)$ such that the spectrum of any $T'\in W_T$ is
disjoint from $I$. 

(b) Let $J\subset \R$ be a finite interval whose boundary points are
no eigenvalues of $T$. Then there exists an open neighbourhood
$W_T$ of $T$ in $\FF_s(E,F)$ such that the assignment
\begin{equation}\label{eq:specprojcont}
  W_T\to\LL(F,F),\qquad T'\mapsto E_J^{T'}
\end{equation}
is continuous.
\end{lemma}

{\bf Proof: }
(a) The family $\{T_t:=T-t\,{\rm Id}\}_{t\in I}\subset \FF_s(E,F)$
consists of invertible operators, and every $T_t$ has an
open neighbourhood $U_t$ which consists of invertible operators. Since 
the family $\{T_t\}_{t\in I}$ is compact, it is covered by finitely many
$U_{t_j}$, $j=1,\dots,l$ and we set
$$
  W_T:=\bigcap_{j=1}^l(U_{t_j}+t_j 
  {\rm Id}).
$$
(b) By part (a), we find an open neighbourhood $W_T$ such that no
$T'\in W_T$ has eigenvalues at the boundary points of $J$. The desired
continuity can now be deduced e.g.~from Dunford calculus~\cite{Dunford} as 
follows. Let $\Gamma\subset \mathbb{C}$ be a smooth simple closed curve
enclosing $J$ and disjoint from the spectrum of $T'$ for all $T'\in W_T$.
We orient $\Gamma$ as the boundary of the connected component of 
$\mathbb{C}\setminus \Gamma$ containing $J$.
Then the spectral measure in question can be expressed by the
operator valued Cauchy integral formula
\begin{equation}\label{eq:dunford}
  E_J^{T'} = \frac{1}{2\pi i}\int_\Gamma(\lambda{\rm Id}-T')^{-1}d\lambda,
\end{equation}
see e.g.~\cite{Dunford}, \S\,X.1 Formula~(i). In loc. cit. it is stated for bounded operators, but it 
immediately generalizes to our case
as follows. Note that the integral itself on the right hand side 
of~\eqref{eq:dunford} still makes sense. Let now $e$ be an eigenvector for an eigenvalue $\mu$
of $T'$. Then 
$$
(\lambda{\rm Id}-T')^{-1}e=
(\lambda-\mu)^{-1}e.
$$
Since the spectrum of $T'$ consists
of eigenvalues, the usual Cauchy formula from complex analysis 
implies that the right hand side 
of~\eqref{eq:dunford} does indeed define the desired spectral measure.
The continuous dependence of $E_J^{T'}$ on $T'$ now follows from 
the continuous dependence of the Cauchy integral on its integrand.
%marginpar{ Check the adjusted argument! EV}
\hfill$\square$

%%%
\subsection{The main result and strategy of proof}\label{ss:main-result}
%%%

We introduce the following notation for any real number 
$\mathfrak{R}$, where $\sigma(T)$ denotes the spectrum of $T$:
\begin{align*}
  \FF_s^{>\mathfrak{R}}(E,F) &:= \{T\in \FF_s(E,F)\mid 
  \sigma(T)\subset (\mathfrak{R},+\infty)\}, \cr
  \FF_s^{>\mathfrak{R}}(E,F)^* &:= \{T\in \FF_s^{>\mathfrak{R}}(E,F) \mid 
  \ker T=\coker T=0 \}, \cr
  \FF_s^{<\mathfrak{R}}(E,F) &:= \{T\in \FF_s(E,F)\mid 
  \sigma(T)\subset (-\infty,\mathfrak{R})\}.
\end{align*}

Now we can state the main result of this section.

\begin{thm}\label{thm:det}
Consider the setting above and assume that the inclusion\\ $\iota:E\into F$
is compact. Let $\mathfrak{R}$ be any real number. 
Then the equalities $(\im T)^{\perp_F} = \ker T$ give rise to a
canonical orientation of the restrictions  
$$
  \det|_{\FF_s^{>\mathfrak{R}}(E,F)}\quad\text{and}\quad
  \det|_{\FF_s^{<\mathfrak{R}}(E,F)}
$$
of the determinant line bundle to the subspaces of index zero
symmetric Fredholm operators bounded from below or above.  
\end{thm}

The proof will occupy the rest of this section. We begin by explaining
the main ideas. Observe that the involution 
$$
\FF_s^{>\mathfrak{R}}(E,F)\longrightarrow 
\FF_s^{< -\mathfrak{R}}(E,F),\quad A\mapsto -A
$$
preserves the determinant line bundle. Therefore, it is enough to
prove our statement for $\FF_s^{>\mathfrak{R}}(E,F)$. 

For $T\in\FF_s(E,F)$ let $v_1,\dots,v_k$ be an
orthonormal basis (with respect to $\la\ ,\ \ra_F$) of 
$(\im T)^{\perp_F} = \ker T$. We identify 
$\ker T$ with its dual using the Hilbert space structure. The element 
\begin{equation}\label{eq:secKai}
   t(T):=v_k\wedge\cdots\wedge v_1\otimes v_1\wedge\cdots\wedge v_k\in
   \Lambda^k(\ker T)\otimes\Lambda^k(\ker T) = \det(T)
\end{equation}
does not depend on the choice of the orthonormal basis, 
so these elements define a canonical section $t$ of
$\det|_{\FF_s(E,F)}$. Unfortunately, this section will turn out to be
discontinuous --- it has singularities at operators that have more
kernel than some of their neighbours. Observe, however, that the restriction of the determinant line bundle to 
the subset $\FF_s(E,F)^*$ 
of invertible operators is a trivial 
bundle $\FF_s(E,F)^*\times \R$. 
Moreover, the section $t$ over 
$\FF_s(E,F)^*$ is the constant 
$1\in \R$. The section $t$ can be thought of as a tautological section 
and its restriction to 
$\FF_s(E,F)^*$ is tautologically continuous.

The idea now is to compensate the above discontinuity using a certain
spectral count.
Recall from Lemma~\ref{lem:spectrum} that the spectrum
$\sigma(T)\subset\R$ of $T\in\FF_s(E,F)$ is discrete and consists only
of eigenvalues of finite multiplicity. Let  
$$
\rho:\R\to (-\infty,1]
$$ 
be a nondecreasing smooth function which restricts 
on the interval $(-\infty,a]$ as the identity
$\lambda\mapsto \lambda$
and which is constant equal to $1$ on 
the interval $[b,\infty)$ for
some $0<a<b$. We define our spectral count
$$
  \mu:\FF_s^{>\mathfrak{R}}(E,F)\longrightarrow \R
$$
by the formula
$$
  \mu(T):=\prod_{\lambda\in \sigma(T)\setminus \{0\}}
\rho(\lambda),
$$
%By convention, if eigenvalue $\lambda$ has multiplicity
%$a$ we list it $a$ times in $\sigma_p(T)$. For example if 
%$\dim\ker T=a$, then we write $\lambda_1=\dots=\lambda_a=0$
%for the first $a$ eigenvalues of $T$.
where each eigenvalue appears with its multiplicity.
Note that the product above contains only finitely many factors
different from $1$ because all eigenvalues of $T$ are $>\mathfrak{R}$
and $\rho(\lambda)=1$ for large $\lambda$. 
The function $\mu$ is nowhere zero. Observe, however, that 
part of the spectrum may converge to $0$ in a family 
of operators. On such a family the function $\mu$ converges 
to $0$ and compensates the singularity in the section $t$ above. 
We define a modified section of 
the determinant line bundle over
$\FF_s^{>\mathfrak{R}}(E,F)$ by
$$
  s(T):=\mu(T)v_k\wedge\cdots\wedge v_1\otimes v_1\wedge\cdots\wedge v_k.
$$
We will show that $s$ is continuous. Note that different 
choices of the function $\rho$ lead to positively 
proportional sections $s$, therefore the induced orientation is canonical.

Let $T\in \FF_s^{>\mathfrak{R}}(E,F)$ be given. We specify the general 
Fredholm setting above to the present situation. Namely, set 
$N:=\dim\ker T$ and pick a linear isomorphism 
$$
  \Phi:\R^N\longrightarrow \ker T=
(\im T)^{\perp F}\subset F
$$
respecting scalar products. Let $U_{T,\Phi}$ be a neighbourhood of $T$
in $\FF(E,F)$ as above and consider its restriction to 
$U_T:=U_{T,\Phi}\cap \FF_s^{>\mathfrak{R}}(E,F)$. 
Let $\{T_n\}_{n\in \N}\subset U_T$ be any sequence with 
$$
  \lim_{n\to\infty}T_n=T.
$$ 
We need to show 
\begin{equation}\label{eq:mainconv1}
  \lim_{n\to\infty}s(T_n)=s(T).
\end{equation}
By a standard trick, it suffices to show that for every
subsequence $(T_{n_k})$ there exists a subsequence $(T_{n_{k_j}})$
such that $\lim_{j\to\infty}s(T_{n_{k_j}})=s(T)$. Using this, in the
proof of~\eqref{eq:mainconv1} we will repeately pass to subsequences
of $(T_n)$, always renaming them back to $T_n$.  

Inequality~\eqref{eq:dimineq} with $k=l$ (because of index zero)
implies that after passing to a subsequence we can assume 
$$
  \dim \ker T_n=k\le \dim \ker T = N
$$
for some constant $k$. We understand the desired 
convergence~\eqref{eq:mainconv1} in terms of the above 
isomorphism $\iota_\Phi$.  That is, we have to prove that 
\begin{equation}\label{eq:keyconv}
  \lim_{n\to\infty}\iota_\Phi\bigl(T_n,s(T_n)\bigr)=
\iota_\Phi\bigl(T,s(T)\bigr).
\end{equation}
Let us recall the choices we have to make in order to 
define the necessary objects entering the last assertion:
\begin{itemize}
 \item [(i)] an orthonormal basis $v_1,\dots,v_N$ of 
 $\ker T$;
 \item [(ii)] an orthonormal basis $v_1^n,\dots,v_k^n$ of 
 $\ker T_n$ for each $n\in \N$;
 \item [(iii)] a basis $\zeta_1,\dots,\zeta_N$ of $\R^N$ 
 (to define $\iota_\Phi(T,\cdot)$);
 \item [(iv)] a basis $\zeta_1^n,\dots,\zeta_N^n$ of $\R^N$ 
 for each $n\in \N$ (to define $\iota_\Phi(T_n,\cdot)$);
 \item [(v)] a collection $\{\xi_j^n\}_{j=k+1}^N\subset E$ subject to 
 \begin{equation}\label{eq:relxin}
 T_n\xi_j^n=-\Phi\zeta_j^n,\quad j=k+1,\dots,N
 \end{equation}
 (to define $\iota_\Phi(T_n,\cdot)$).
 \end{itemize}
This allows us to write out 
$$
\iota_\Phi\bigl(T_n,s(T_n)\bigr)= (-1)^kC_n\mu(T_n)
(v_1^n,0)\wedge\cdots\wedge(v_k^n,0)\wedge 
(\xi_{k+1}^n,\zeta_{k+1}^n)\wedge\cdots\wedge
(\xi_N^n,\zeta_N^n)
$$
and
$$
\iota_\Phi\bigl(T,s(T)\bigr)= (-1)^NC\mu(T)
(v_1,0)\wedge\cdots\wedge(v_k,0)\wedge 
(v_{k+1},0)\wedge\cdots\wedge
(v_N,0)
$$
with the abbreviations
$$
  C_n:=\frac{\det\la v_i^n,\Phi\zeta_j^n\ra
_{i,j=1,\dots,k}}{\det(\zeta_1^n,\dots,\zeta_N^n)},\qquad
  C:=\frac{\det\la v_i,\Phi\zeta_j\ra
_{i,j=1,\dots,N}}{\det(\zeta_1,\dots,\zeta_N)}.
$$
The desired convergence~\eqref{eq:keyconv}
thus follows from the following three assertions (for suitable
choices (i)--(v)):
%\begin{itemize}
% \item [(A)] 

Assertion {\bf(A)}
 $$
 \lim_{n\to\infty}v_1^n\wedge\cdots\wedge v_k^n=
 v_1\wedge\cdots\wedge v_k.
 $$ 
% \item [(B)]

Assertion {\bf(B)}
 $$\lim_{n\to\infty}\mu(T_n)(-\xi_{k+1}^n,\zeta_{k+1}^n)\wedge\cdots\wedge
 (-\xi_N^n,\zeta_N^n)=\mu(T)(v_{k+1},0)
 \wedge\cdots\wedge(v_N,0).
 $$
% \item [(C)]

Assertion {\bf(C)}
\begin{equation}\label{eq:detconv1}
 \lim_{n\to\infty}\det\la v_i^n,\Phi\zeta_j^n\ra_{i,j=1,\dots,k}= 
 \det\la v_i,\Phi\zeta_j\ra_{i,j=1,\dots,N},
 \end{equation}
\begin{equation}\label{eq:detconv2}
 \lim_{n\to\infty}\det(\zeta_1^n,\dots,\zeta_N^n)=\det(\zeta_1,\dots,\zeta_N).
 \end{equation}
%\end{itemize} 
Let us describe the strategy for making the choices above.
We are going to show that in a certain sense $\ker T_n$
converges to a subspace 
$\ker_{\ker}T$ of $\ker T$ as 
$n\to\infty$. Note that each $\ker T_n$ is a member of an infinite
dimensional Grassmannian. Therefore, we need an additional construction 
(``parametrizing the kernels'') relating $\ker T_n$
to a subspace of $\ker T$. This will give us the key orthogonal splitting 
\begin{equation*}%\label{eq:keysplitintro}
  \ker T=\ker_{\ker}T\oplus \ker_{\im}T,
\end{equation*}
where $\ker_{\im}T$ is the orthogonal complement 
to $\ker_{\ker}T$ in $\ker T$. Alternatively, $\ker_{\im}T$
is the limit of $\im T_n\cap \ker T$ as $n\to\infty$.
We let $v_1,\dots,v_k$
be any orthonormal basis of $\ker_{\ker}T$ and 
$v_{k+1},\dots,v_N$ be any orthonormal basis of 
$\ker_{\im}T$. We define $\{\zeta_j\}_{j=1}^N$ to be 
the preimage of $\{v_j\}_{j=1}^N$ under $\Phi$
and $\zeta_j^n:=\zeta^j$ for $j=1,\dots,k$.
We then define $\{v_j^n\}_{j=1}^k$ as the orthonormalization of a
suitable projection of $\{v_j\}_{j=1}^k$ to $\ker T_n\cap \ker T$, 
and $\{\Phi\zeta_j^n\}_{j=k+1}^N$ as the orthogonal projection of 
$\{v_j\}_{j=k+1}^N$ to $\im T_n$. This will allow us to define
$\{\xi_j\}_{j=k+1}^N$ using~\eqref{eq:relxin}
%, where $T_n$ is viewed as an isomorphism between $\im T_n\cap E$ and
%$\im T_n$. 
%In order to work with the above convergences comfortably we need some
%statements about the spectral measures of the operators
%involved. Therefore, the first step is to show that symmetric index
%zero Fredholm operators are in fact self-adjoint. Once this is
%established and the kernels of $T_n$ are suitably parametrized, we
%will begin the actual proof of the three assertions above.
and it remains to verify the three assertions above (see Section~\ref{ss:proofthmdet}).

%%%
\subsection{Parametrizing the kernels}
%%%

Let $X$ and $Y$ be two Banach spaces and $\FF^*(X,Y)$
be the space of surjective Fredholm 
operators from $X$ to $Y$.
Pick any $D^0\in \FF^*(X,Y)$ and let 
$U_{D^0}\subset \FF^*(X,Y)$ be a connected open 
neighbourhood of $D^0$. We describe one possible trivialization of the
kernel bundle over $U_{D^0}$. 
Consider the direct sum $Y\oplus \ker D^0$ and let 
$\C_0\subset X$ be a closed direct complement of $\ker D^0$ in $X$. Let 
$$
  P_{\ker D^0}:X\longrightarrow \ker D^0
$$
denote the projection onto $\ker D^0$ along $\C_0$.

To every operator  
$D\in U_{D^0}$ we associate an operator 
$\tilde D\in L(X,Y\oplus \ker D^0)$ defined by 
$\tilde D(x):=(D(x),P_{\ker D^0}(x))$. By construction the operator
$\tilde D^0$ is bijective. Therefore, by shrinking  
$U_{D^0}$ if necessary, we may assume that $\tilde D$ is bijective for
every $D\in U_{D^0}$. As a consequence, each $\tilde D$ 
admits continuous inverse depending continuously on $D$.
%This allows us to construct a family of isomorphisms 
%$$
%\QQ_{D}:\ker D^0\longrightarrow \ker D
%$$
%depending continuously on $D\in U_{D^0}$ and such that 
%$\QQ_{D^0}={\rm Id}_{\ker D^0}$. Indeed,
Given any $x\in X$, let $\tilde x\in X$ be the unique solution of the
equation  
\begin{equation}\label{eq:tilde x}
\tilde D(\tilde x)=(D(x),0).
\end{equation}
Writing out equation~\eqref{eq:tilde x} in components gives us 
$$
  D(\tilde x)=D(x),\quad P_{\ker D^0}(\tilde x)=0.
$$
In other words,
$$
  D(x-\tilde x)=0,\quad \tilde x\in \C_0.
$$
Observe that $\tilde x$ depends continuously on 
$(D,x)\in U_{D^0}\times X$, and $\tilde x=0$
for $D=D^0$ and any $x\in\ker D^0$. So the isomorphisms
$$
  \QQ_D:\ker D^0\stackrel{\cong}\longrightarrow \ker D,\qquad \QQ_D(x):=x-\tilde x
$$
depend continuously on $D\in U_{D^0}$ and 
$$
%||\QQ_D-id||\to 0\quad \text{as}\quad D\to D^0.  
  \QQ_{D^0}={\rm Id}_{\ker D^0}. 
$$
To apply this in our situation 
set $X:=E\oplus \R^N$, $Y:=F$, $D^0:=T\oplus \Phi$, 
$D:=T_n\oplus \Phi$. Furthermore, in our case 
$\ker T\cong \ker T\times \{0\}=\ker (T\oplus \Phi)$
and we can take $\C_0:=\im T\cap E$. Then the preceding construction
gives us isomorphisms
$$
  \QQ_n:=\QQ_{T_n\oplus \Phi}:\ker T\stackrel{\cong}\longrightarrow 
\ker (T_n\oplus \Phi)
$$
such that
\begin{equation}\label{eq:Qconv}
  \lim_{n\to\infty}\QQ_n = {\rm Id}_{\ker T}.
\end{equation}
Equation~\eqref{eq:Qconv} implies in particular 
\begin{equation}\label{eq:nconv}
  \|(\QQ_n)_1-{\rm Id}_{\ker T}\|\to 0,\qquad 
  \|\QQ_n^{-1}|_{\ker T_n\times\{0\}} - {\rm Id}_{\ker T_n\times\{0\}}\|\to 0
\end{equation}
as $n\to\infty$, where $(\QQ_n)_1$ denotes the first component of $\QQ_n$.

Using the isomorphisms $\QQ_n$ we can now establish the orthogonal
splitting $\ker T=\ker_{\ker} T\oplus\ker_{\im} T$.  
%The following quantity measures the ``angle'' between two subspaces of
%a finite dimensional space.
For two subspaces $H_1,H_2$ of a finite dimensional Hilbert space
$(H,\la\cdot,\cdot\ra_H)$ we define
$$
  \la H_1,H_2 \ra := \sup\{\la u,v\ra_H\mid (u,v)\in H_1\times
  H_2,\,\|u\|=\|v\|=1\}.
$$
Note that $\la H_1,H_2 \ra=0$ if and only if $H_1$ and 
$H_2$ are orthogonal to each other.
We consider $\ker T$ and $\ker(T_n\oplus\Phi)$ as finite dimensional
Hilbert spaces with the scalar products induced from $F$. 
Let ``$\perp$'' denote ``orthogonal complement in
$\ker(T_n\oplus\Phi)$'' for the discussion below. For a subspace $K$
of $\ker(T_n\oplus\Phi)$ let $(K)_2\subset \R^N$ denote the projection
of $K$ on the second component. Observe that  
$$
  (\ker T_n\times \{0\})^\perp=\{(\xi,\zeta)\in 
 (\im T_n\cap E)\times \R^N\mid T_n\xi=-\Phi\zeta\}.
$$
Since $T_n$ restricts to an isomorphism from $\im T_n\cap E$
to $\im T_n$, we get 
$$
  \Phi((\ker T_n\times \{0\})^\perp)_2 = \im T_n\cap \ker T.
$$
Since the Grassmannian of $k$-dimensional subspaces of $\ker T$ is compact, we get
(after passing to a subsequence if necessary) the following limits: 
\begin{equation}\label{eq:defUV}
  \ker_{\ker} T:=\lim_{n\to \infty}Q_n^{-1}(\ker T_n\times \{0\}),\qquad
  \ker_{\im} T:=\lim_{n\to \infty}(\im T_n\cap \ker T). 
\end{equation}
For a closed subspace $C\subset F$ we denote by
$$
  P_C:F\to F
$$
the orthogonal projection onto $C$. In terms of orthogonal
projections, the preceding equations mean
\begin{equation}\label{eq:defUVproj}
  P_{\ker_{\ker} T} = \lim_{n\to \infty}P_{Q_n^{-1}(\ker T_n\times \{0\})},\qquad
  P_{\ker_{\im} T} = \lim_{n\to \infty}P_{\im T_n\cap \ker T}. 
\end{equation}
Our next goal is to show orthogonality of $\ker_{\ker} T$
and $\ker_{\im} T$.
Note that $\ker T_n\times \{0\}$ and 
$\im T_n\cap \ker T$ are orthogonal to each other as subspaces of
$F$. This together with equation~\eqref{eq:nconv} implies that  
$$
  \la\im T_n\cap \ker T,Q_n^{-1}(\ker T_n\times \{0\})\ra\to 0\quad
  \text{as}\quad n\to\infty. 
$$
The desired orthogonality now follows from equation~\eqref{eq:defUV}. Since 
$$
  \dim\ker_{\ker} T=k\quad\text{and}\quad \dim\ker_{\im} T=N-k,
$$
we get the orthogonal splitting
\begin{equation}\label{eq:keysplitintro}
\ker T=\ker_{\ker} T\oplus\ker_{\im} T.
\end{equation}

%%%
\subsection{Proof of Theorem~\ref{thm:det}}\label{ss:proofthmdet}
%%%

{\bf The choices (i)--(v). }
Now we are ready to make the choices (i)--(v) in Section~\ref{ss:main-result}.
For (i) let $v_1,\dots,v_N$ be an orthonormal basis of 
$\ker T$ such that $v_1,\dots,v_k$ form a basis of 
$\ker_{\ker}T$, and $v_{k+1},\dots,v_N$ form a basis of
$\ker_{\im}T$. For (ii), let $GS_n$ denote the Gram-Schmidt retraction from 
the space of all bases of $\ker T_n$ to the space of orthonormal bases
of $\ker T_n$ and set
\begin{equation}\label{eq:defbasisvn}
  \{v_j^n\}_{j=1}^k := GS_n\Bigl(P_{\ker T_n\times \{0\}}\QQ_n(\{v_j\}_{j=1}^k)\Bigr).
\end{equation}
For (iii) we define 
\begin{equation}\label{eq:defzeta}
\zeta_j:=\Phi^{-1}v_j,\quad j=1,\dots,N.
\end{equation}
For (iv) we define 
\begin{equation}\label{eq:defzetan}
  \zeta_j^n := \begin{cases}
    \zeta_j & 1\leq j\leq k, \cr
    \Phi^{-1}P_{\im T_n\cap \ker T}v_j & k+1\leq j\leq N.
  \end{cases}
\end{equation}
The second part of this definition rewrites as
\begin{equation}\label{eq:zetan}
  \Phi\zeta_j^n=P_{\im T_n\cap \ker T}v_j,\quad j=k+1,\dots,N.
\end{equation}
In view of the second equation in~\eqref{eq:defUVproj} and $v_j\in
\ker_{\rm im}T$ for $j=k+1,\dots,N$ this implies 
\begin{equation}\label{eq:zetanconv}
  \Phi\zeta_j^n\to P_{\ker_{\rm im}T}v_j = v_j\quad\text{as
  }n\to\infty,\qquad j=k+1,\dots, N. 
\end{equation}
%We are going to extract two consequences of the last displayed equation. For the first one
Applying $\Phi^{-1}$ to both sides of the last equation, we get 
\begin{equation}\label{eq:zetaconv}
\zeta_j^n\to \zeta_j\quad\text{as}\quad n\to\infty,\quad j=k+1,\dots N.
\end{equation}
In particular, we see that $\{\zeta_j^n\}_{j=1}^N$ form a 
basis of $\R^N$ for large $n$.
For (v), we use the isomorphism
\begin{equation}\label{eq:Tntilde}
  \wt T_n := T_n|_{\im T_n\cap E}:\im T_n\cap E\stackrel{\cong}\longrightarrow \im T_n  
\end{equation}
to define 
\begin{equation}\label{eq:xijn}
  \xi_j^n := -\wt T_n^{-1}\Phi\zeta_j^n,\quad j=1,\dots,N.
\end{equation}
%The second consequence of~\eqref{eq:zetanconv} will be stated later
%in the proof of part (B). 

With these choices, we will now verify assertions (A)--(C) from
Section~\ref{ss:main-result}.  
\smallskip

{\bf Proof of Assertion (A). }
We will use the following simple lemma. 

\begin{lemma}\label{lem:convgen}
Let $\{A_n\}_{n\in \N}$ and $\{B_n\}_{n\in \N}$
be two bounded sequences of operators $A_n\in\LL(Y,Z)$,
$B_n\in\LL(X,Y)$ between Banach spaces. 
Assume convergence $B:=\lim_{n\to\infty}B_n $
and the equality $A_nB_n=0$ for all $n\in \N$. Then 
$\lim_{n\to\infty}A_nB=0$. 
\end{lemma}

{\bf Proof:}
Since the sequence $\{A_n\}_{n\in \N}$ is bounded, we get 
$$
  \|A_m(B-B_n)\|\to 0\quad \text{as} \quad n\to\infty
$$
uniformly with respect to $m$. Set $m:=n$ and use $A_nB_n=0$
to conclude.
\hfill $\square$

We apply the lemma with $X=Y=Z=\ker T$ and the projections
$$
  A_n := \QQ_n^{-1}P_{(\ker T_n\times \{0\})^\perp}\QQ_n,\qquad 
  B_n := P_{\QQ_n^{-1}(\ker T_n\times \{0\})}.
$$ 
Since obviously $A_nB_n=0$, and $\lim_{n\to\infty}B_n=P_{\ker_{\ker} T}$ by
  equation~\eqref{eq:defUVproj}, Lemma~\ref{lem:convgen} yields
\begin{equation}\label{eq:convUn}
\lim_{n\to\infty} \QQ_n^{-1}P_{\ker T_n\times \{0\}^\perp}\QQ_nP_{\ker_{\ker}T}=0.
\end{equation}
In view of $P_{\ker T_n\times \{0\}^\perp}+P_{\ker T_n\times \{0\}}={\rm
  Id}_{\ker T}$ and equation~\eqref{eq:Qconv}, this implies
\begin{equation}\label{eq:convVn}
\lim_{n\to\infty} P_{\ker T_n\times \{0\}}\QQ_nP_{\ker_{\ker}T}=P_{\ker_{\ker}T}.
\end{equation}
This implies
\begin{equation}\label{eq:convprebasisvn}
  P_{\ker T_n\times \{0\}}\QQ_n(\{v_j\}_{j=1}^k)\to
  \{v_j\}_{j=1}^k\quad \text{as}\quad n\to\infty, 
\end{equation}
and since $\{v_j\}_{j=1}^k$ is orthonormal we get 
\begin{equation}\label{eq:convbasisvn}
  \{v_j^n\}_{j=1}^k\to \{v_j\}_{j=1}^k \quad \text{as}\quad n\to\infty.
\end{equation}
This proves Assertion (A). 
\hfill $\square$
\smallskip

{\bf Proof of Assertion (C). }
The choice~\eqref{eq:defzeta} of $\{\zeta_j\}_{j=1,\dots,N}$ implies that 
the right hand side of~\eqref{eq:detconv1} equals $1$. 
The choice~\eqref{eq:defzetan} of $\{\zeta_j^n\}_{j=1,\dots,k}$ and 
equation~\eqref{eq:convbasisvn} implies that the left hand side
of~\eqref{eq:detconv1} equals $1$. 
Equation~\eqref{eq:detconv2} follows from
equations~\eqref{eq:defzetan} (for $j\leq k$) and~\eqref{eq:zetaconv}
(for $j>k$).
\hfill $\square$
\smallskip

{\bf Proof of Assertion (B). }
First, we prepare several convergence statements.
Let $\eps>0$ be such that $0$ is the only 
eigenvalue of $T$ in $(-2\eps,2\eps)$
and the function $\rho$ equals the identity on $(-\eps,\eps)$.
Then by continuity~\eqref{eq:specprojcont} of the spectral measure we get
\begin{equation}\label{eq:specconte}
  \lim_{n\to\infty}E_{(-\eps,\eps)}^{T_n}=
E_{(-\eps,\eps)}^T=P_{\ker T},\qquad 
\lim_{n\to\infty}E_{\R\setminus(-\eps,\eps)}^{T_n}
%=E_{\R\setminus(-\eps,\eps)}^T
={\rm Id}-P_{\ker T}.
\end{equation}
Let $k+1\leq j\leq N$. We apply the operations in~\eqref{eq:specconte}
to both sides of equation~\eqref{eq:zetanconv} to get 
\begin{equation}\label{eq:phizetanconv}
  \lim_{n\to\infty}E_{(-\eps,\eps)}^{T_n}\Phi\zeta_j^n=v_j,\qquad
  \lim_{n\to\infty}E_{\R\setminus(-\eps,\eps)}^{T_n}\Phi\zeta_j^n=0.
\end{equation}
Observe that $P_{\ker T_n}\Phi\zeta_j^n=0$ and define
$$
  z_j^n:=E_{(-\eps,\eps)\setminus \{0\}}^{T_n}
\Phi\zeta_j^n=E_{(-\eps,\eps)}^{T_n}\zeta_j^n,\qquad 
\tilde z_j^n:=
E_{\R\setminus(-\eps,\eps)}^{T_n}\Phi\zeta_j^n.
$$
Note that
$$
\Phi\zeta_j^n=z_j^n+\tilde z_j^n.
$$
Recall the isomorphism $\wt T_n$ from~\eqref{eq:Tntilde}. 
The absolute values of the eigenvalues of 
$\wt T_{n,\eps}^{-1}:=\wt T_n^{-1}|_{\im E^{T_n}_{\R\setminus(-\eps,\eps)}}$ 
are bounded above by $\eps^{-1}$. Since the spectral radius of 
%$T_n^{-1}|_{\im E^{T_n}_{\R\setminus(-\eps,\eps)}}$ 
a compact self-adjoint operator equals its norm, this implies
\begin{equation}\label{eq:restest}
  \|\wt T_n^{-1}\tilde z_j^n\| = \|\wt T_{n,\eps}^{-1}\tilde z_j^n\|\le 
  \|\wt T_{n,\eps}^{-1}\|\,\|\tilde z_j^n\|\le \eps^{-1}\|\tilde z_j^n\|.
\end{equation}
The preceding discussion can be summarized as follows for $j=k+1,\dots,N$:
\begin{equation}\label{eq:summary}
  \Phi\zeta_j^n=z_j^n+\tilde z_j^n,\qquad
  \lim_{n\to\infty}z_j^n=v_j,\qquad
  \lim_{n\to\infty}\tilde z_j^n=\lim_{n\to\infty}\wt T_n^{-1}\tilde z_j^n=0.
\end{equation}
Second, we manipulate the spectral count $\mu$ and 
determinants of certain finite dimensional operators.
%$$
%\det (T_n^{-1}|_{\im E^{T_n}_{(-\eps,+\eps)\setminus \{0\}}})
%=\prod_{\lambda\in \sigma(T_n)\cap(-\eps,\eps)\setminus\{0\}}\lambda^{-1}
%$$
Set 
$$
  \mu_n:=\prod_{\lambda\in\sigma(T_n)\cap(\R\setminus(-\eps,\eps))}\rho(\lambda).
$$
Since $\rho(\lambda)=\lambda$ for $\lambda\in(-\eps,\eps)$, we get 
\begin{equation}\label{eq:detmu}
  \mu(T_n)=\mu_n \prod_{\lambda\in \sigma(T_n)\cap((-\eps,\eps)\setminus\{0\})}\lambda
  = \mu_n\det (\wt T_n|_{\im E^{T_n}_{(-\eps,+\eps)\setminus \{0\}}})
\end{equation}
and 
\begin{equation}\label{eq:convmun}
\lim_{n\to\infty}\mu_n=\mu(T).
\end{equation}
%For future use we reserve
%\begin{equation}\label{eq:detmu}
%\mu(T_n)\det(T_n^{-1}|_{\im E^{T_n}_{(-\eps,+\eps)\setminus \{0\}}})=\mu_n. 
%\end{equation}
Now we attend to the main part of the argument.
Using the notation $\Lambda_{j=a}^bw_j=w_a\wedge\cdots\wedge w_b$,
Assertion (B) reads
$$
  \lim_{n\to\infty}\mu(T_n)\Lambda_{j=k+1}^N(-\xi_j^n,\zeta_j^n) = \mu(T)\Lambda_{j=k+1}^N(v_j,0).
$$
Identifying vectors $\xi\in E$ and $\zeta\in\R^N$ with their images
$(\xi,0)$ and $(0,\zeta)$ in $E\oplus\R^N$ and using
equations~\eqref{eq:xijn} and~\eqref{eq:summary}, this reads
\begin{equation}\label{eq:B}
  \lim_{n\to\infty}\mu(T_n)\Lambda_{j=k+1}^N(\wt T_n^{-1}z_j^n + \wt
  T_n^{-1}\tilde z_j^n + \zeta_j^n) = \mu(T)\Lambda_{j=k+1}^Nv_j.
\end{equation}
We split the wedge product on the left hand side as the leading term plus the rest, 
\begin{equation}\label{eq:split}
   \Lambda_{j=k+1}^N(\wt T_n^{-1}z_j^n + \wt
  T_n^{-1}\tilde z_j^n + \zeta_j^n) = \Lambda_{j=k+1}^N\wt T_n^{-1}z_j^n + \tilde R_n.
\end{equation}
%where 
%$$
%\tilde R_n:=\sum\Lambda_{j=k+1}^NT_n^{-1}w_j^n,
%$$
%and the sum is taken over all choices $w_j^n\in\{z_j^n,\tilde
%z_j^n\}$ such that for at least one $j\in \{k+1,\dots,N\}$ we have
%$w_j^n=\tilde z_j^n$.
We first discuss the rest term 
$\tilde R_n$. Modulo signs a typical summand is
%$$
%(T_n^{-1}z_{j_{k+1}}^n\wedge\cdots\wedge T_n^{-1}z_{j_l}^n)\wedge
%(T_n^{-1}\tilde z_{j_{l+1}}^n\wedge\cdots\wedge T_n^{-1}\tilde z_{j_N}^n)
%$$
$$
   \Bigl(\Lambda_{i=k+1}^l\wt T_n^{-1}z_{j_i}^n\Bigr)\wedge
   \Bigl(\Lambda_{i=l+1}^m\wt T_n^{-1}\tilde z_{j_i}^n\Bigr)\wedge
   \Bigl(\Lambda_{i=m+1}^N\zeta_{j_i}^n\Bigr)
$$
for some $k\leq l\leq m\leq N$ with $l<N$. 
Equation~\eqref{eq:zetanconv} and the third equation
in~\eqref{eq:summary} imply that the second and third factors in the last
displayed equation
%converges to $0$
remain bounded as $n\to\infty$. For the first factor, we pick 
orthonormal bases $e^n_{k+1},\dots,e^n_N$ of the spaces
$\im E^{T_n}_{(-\eps,+\eps)\setminus \{0\}}$ consisting of 
eigenvectors of $T_n$ and converging to an orthonormal system
$e_{k+1},\dots,e_N$ as $n\to\infty$. 
We write $z_j^n=\sum_{i=k+1}^Nc_{jn}^ie_i^n$ in these bases, with 
coefficients $c_{jn}^i\in \R$. Since $z_j^n\to v_j$ by
the second part of~\eqref{eq:summary},
%implies boundedness of the sequence $\{z_j^n\}_{n\in \N}$ for all $j\in \{k+1,N\}$.
the sequence $\{c_{jn}^i\}_{n\in \N}$ converges and thus
in particular remains bounded as $n\to\infty$, for all
$i,j=k+1,\dots,N$. This way the first factor
$\wt T_n^{-1}z_{j_1}^n\wedge\cdots\wedge\wt T_n^{-1}z_{j_l}^n$
can be viewed as a homogeneous polynomial in 
$\lambda_{k+1}^{-1},\dots,\lambda_N^{-1}$
of degree $l-k$ strictly less than $N-k$. 
Moreover, this polynomial has bounded coefficients in 
$\Lambda^{l-k}\im E_{(-\eps,+\eps)\setminus\{0\}}^{T_n}$,
and each variable $\lambda_j^{-1}$ for $j=k+1,\dots,N$
enters every monomial with power $0$ or $1$. This implies 
\begin{equation}\label{eq:partconv}
  \lim_{n\to\infty}\mu(T_n)\Bigl(\Lambda_{i=k+1}^l\wt T_n^{-1}z_{j_i}^n\Bigr)=0,
\end{equation}
and therefore
\begin{equation}\label{eq:restconv}
  \lim_{n\to\infty}\mu(T_n)\tilde R_n=0.
\end{equation}
Now we consider the leading term in equation~\eqref{eq:split} which we
rewrite as
\begin{align*}%\label{eq:keylead}
  \mu(T_n)\Lambda_{j=k+1}^N\wt T_n^{-1}z_j^n = 
  \mu(T_n)\det(\wt T_n^{-1}|_{\im E^{T_n}_{(-\eps,+\eps)\setminus \{0\}}})\Lambda_{j=k+1}^Nz_j^n =
  \mu_n\Lambda_{j=k+1}^Nz_j^n.
\end{align*}
where the first equality follows from the definition of the
determinant of a finite dimensional operator and the second 
equality follows from equation~\eqref{eq:detmu}.
In view of equation~\eqref{eq:convmun} and the second equation
in~\eqref{eq:summary}, this implies
\begin{equation*}%\label{eq:leadconv}
  \lim_{n\to\infty}\mu(T_n)\Lambda_{j=k+1}^NT_n^{-1}z_j^n=
  \mu(T)\Lambda_{j=k+1}^Nv_j,
\end{equation*}
which together with~\eqref{eq:restconv} proves Assertion (B) and thus
Theorem~\ref{thm:det}. 
%\begin{equation}\label{eq:1stcompconv}
%\lim_{n\to\infty}\mu(T_n)\Lambda_{j=k+1}^N(-\xi_j^n)=
%\mu(T)\Lambda_{j=k+1}^Nv_j.
%\end{equation}
%This together with boundedness of $\zeta_j^n$, $j=k+1,\dots, N$ implies (B).
\hfill $\square$

\begin{rem}\label{rem:invertible-or}
%{\bf Explicit computation for invertible operators.}
Recall that on the subset $\FF_s^{>\mathfrak{R}}(E,F)^*$  
consisting of invertible operators we have another orientation of the
determinant bundle given by the section $t$, see~\eqref{eq:secKai} and
the discussion following it. Observe that for
$T\in\FF_s^{>\mathfrak{R}}(E,F)^*$ the values $s(T)$ and  
$t(T)$ differ by $\mu(T)$, where $\mu(T)$ is the product of all
negative eigenvalues of $T$ and some positive factor resulting from
the restriction $\rho|_{(0,+\infty)}$. 
We rephrase this in terms of the respective orientations. 
Let us call the orientation of $\det|_{\FF_s^{>\mathfrak{R}}(E,F)^*}$
induced by $t$ the {\em tautological orientation}. Any other orientation 
of $\det|_{\FF_s^{>\mathfrak{R}}(E,F)^*}$ differs from the
tautological one by a sign. For $T\in\FF_s^{>\mathfrak{R}}(E,F)$ let
$i(T)$ denote the (finite) number of negative eigenvalues of $T$.
Then the restriction of the orientation $s$ to
$\FF_s^{>\mathfrak{R}}(E,F)^*$ is related to the tautological
orientation by 
$$
  s(T) = (-1)^{i(T)}t(T)\quad\text{for } T\in\FF_s^{>\mathfrak{R}}(E,F)^*.
$$
\end{rem}

%%%
\subsection{A counterexample}
%%%

In this subsection we construct a loop of self-adjoint Fredholm
operators with unbounded spectrum and nonorientable determinant
bundle, thus showing that the boundedness of the spectrum from below
or above in Theorem~\ref{thm:det} is essential.

We consider 
$$
  l^2:=\{x=(x_n)_{n\in \Z}\mid x_n\in\R,\;\sum_{n\in \Z}x_n^2<\infty\}
$$
as a real Hilbert space with scalar product
$$
  \la x,y\ra:=\sum_{n\in \Z}x_ny_n.
$$
Similarly, we consider 
$$
  h^1:=\{x=(x_n)_{n\in \Z}\mid x_n\in\R,\;\sum_{n\in \Z}n^2x_n^2<\infty\}
$$
as a real Hilbert space with scalar product
$$
  \la x,y\ra_{h^1}:=\sum_{n\in \Z}x_ny_n+\sum_{n\in \Z}n^2x_ny_n.
$$
Since $h^1\subset l^2$ is dense and the inclusion $h^1\into l^2$ is
compact, the pair $(h^1,l^2)$ satisfies the hypotheses of
Theorem~\ref{thm:det}. Our goal is to prove

\begin{prop}\label{prop:counterex}
The restriction $\det|_{\FF_s(h^1,l^2)}$ of the determinant line
bundle to the space of self-adjoint Fredholm operators $l^2\supset
h^1\to l^2$ is non-orientable. 
\end{prop}

The proof uses the following lemma. Recall the standard Hilbert space
basis $\{e_n:=(\dots 0,1,0,\dots)\}_{n\in \Z}$ of
$l^2$, where the $1$ occupies position $n$.

\begin{lemma}\label{lem:Bernd}
There exists a continuous path $\{U_\tau\}_{\tau\in [1,2]}$
of unitary isomorphisms of $l^2$ that restrict to (non-unitary)
isomorphisms of $h^1$ such that  
\begin{itemize}
 \item [(i)] $U_1$ maps $e_n$ to $e_{n-1}$ for all $n\in \Z$;
\item [(ii)] $U_2={\rm Id}$;
\item [(iii)] $U_\tau^{-1}(e_0)$ is a positive linear combination of
  $e_0$ and $e_1$ for all $\tau\in (1,2)$.
\end{itemize}
\end{lemma}

{\bf Proof: }
The following proof was found by Bernd Schmidt. 
The idea is to concatenate the following two rotational isotopies
defined in terms of the standard basis of $l^2$. For 
$\tau\in [0,1/2]$ we define a unitary operator $V_\tau^1$ of $l^2$ by
requiring that it sends
\begin{align*}
\begin{cases}
e_{-n}&\mapsto \cos(\pi \tau)e_{-n}-\sin(\pi \tau)e_{n+1}, \\
e_{n+1}&\mapsto \sin(\pi \tau)e_{-n}+\cos(\pi \tau)e_{n+1}
\end{cases}
\end{align*}
for all $n\ge 0$.
Similarly, for $\tau\in [0,1/2]$ we define a unitary operator $V_\tau^2$
of $l^2$ by requiring that it sends
\begin{align*}
\begin{cases}
e_0&\mapsto e_0,\\
e_{-n}&\mapsto \cos(\pi \tau)e_{-n}+\sin(\pi \tau)e_n,\\
e_n&\mapsto -\sin(\pi \tau)e_{-n}+\cos(\pi \tau)e_n
\end{cases}
\end{align*}
for all $n>0$. It is clear that both operators $V_\tau^1$
and $V_\tau^2$ restrict to isomorphisms of $h^1$.
Define 
$$
V_\tau := \begin{cases}
  V_\tau^1 & \tau\in [0,1/2],\\
  V_{\tau-1/2}^2\circ V_{1/2}^1 & \tau\in [1/2,1].
\end{cases}
$$
For $\tau=0$ we have $V_0={\rm Id}$. For $\tau=1$ and $n\geq 0$: 
\begin{align*}
  V_1e_{-n} &= V_{1/2}^2V_{1/2}^1e_{-n} = -V_{1/2}^2e_{n+1} = e_{-(n+1)},\cr
  V_1e_{n+1} &= V_{1/2}^2V_{1/2}^1e_{n+1} = V_{1/2}^2e_{-n} = e_{n}.
\end{align*}
To prove the statement about the preimage of $e_0$, we compute for $\tau\in[0,1/2]$:
\begin{align*}\label{eq:ker2part*}
  &V_\tau(\cos(\pi \tau) e_0+ \sin(\pi \tau) e_1) 
  = V_\tau^1(\cos(\pi \tau) e_0+ \sin(\pi \tau) e_1) \cr 
  &= \cos(\pi \tau) (\cos(\pi \tau)e_0-\sin(\pi \tau)e_1)+
  \sin(\pi \tau) (\sin(\pi \tau)e_0+\cos(\pi \tau)e_1)
  = e_0.
\end{align*}
Since $V_\tau^2$ preserves $e_0$ for all $\tau\in [0,1/2]$, it follows
that $V_\tau^{-1}(e_0)$ is a positive linear combination of
$e_0$ and $e_1$ for all $\tau\in (0,1)$. Therefore, backward
reparametrized family $U_\tau:=V_{2-\tau}$ for $\tau\in [1,2]$ has the
desired properties. 
\hfill $\square$
\smallskip

{\bf Proof of Proposition~\ref{prop:counterex}: }
To construct our counterexample, we consider $L^2([0,1],\mathbb{C})$ 
as a real Hilbert space with scalar product 
$$
  \la f,g\ra := Re\int_0^1f(t)\bar g(t)dt.
$$
Similarly, we consider $H^1([0,1],\mathbb{C})$ 
as a real Hilbert space with scalar product
$$
  \la f,g\ra_{H^1} := Re\int_0^1f(t)\bar g(t)dt+
Re\int_0^1f'(t)\bar g'(t)dt.
$$
We abbreviate from now on
$$
L^2:=L^2([0,1],\mathbb{C})\quad\text{and}\quad  
H^1:=H^1([0,1],\mathbb{C}).
$$
Consider the following densely defined operator in $L^2$: 
$$
  L := -i\frac{d}{dt}:L^2\supset H^1\longrightarrow L^2
$$
and the following loop of dense subspaces of $L^2$:
$$
E_\tau:=\{f\in H^1\mid f(0)\in \R,\,\, e^{i\pi \tau}f(1)\in 
\R\},\quad \tau\in [0,1].
$$
Note that $E_0=H^1$ and the inclusion $E_\tau\into L^2$ is compact for
all $\tau\in[0,1]$. We define the restrictions
$$
  L_\tau := L|_{E_\tau}: L^2\supset E_\tau\to L^2.
$$
We claim that $L_\tau$ is symmetric. Indeed, 
$$
  \la L_\tau f,g\ra=Re\int_0^1-if'(t)\bar g(t)dt
  =Re\left(-if\bar g\bigg|_0^1\right)-Re\int_0^1-if'(t)\bar g(t)'dt.
$$
Now $f(0)\bar g(0)\in \R$ and $f(1)\bar g(1)=
e^{i\pi\tau}f(1)\overline{e^{i\pi\tau}g(1)}\in \R$ implies
$Re\left(-if\bar g\big|_0^1\right)=0$, which in view of 
$if\bar g'=-f\overline{ig'}$ yields
$$
  \la L_\tau f,g\ra=\la f,L_\tau g\ra.
$$
The spectrum of $L_\tau$ consists of eigenvalues
$$
   \lambda_n^\tau = \pi(n-\tau),\qquad n\in\Z
$$
with $1$-dimensional eigenspaces spanned by the eigenvectors
$$
   e_n^\tau(t) = e^{i\pi(n-\tau)t},\qquad n\in\Z. 
$$
In particular, the spectrum of $L_\tau$ is unbounded from both sides.
Note that $\ker L_\tau=0$ for $\tau\in(0,1)$ and $\ker L_0=\R e_0^0$.
It is easy to see that $(\im L_\tau)^\perp=\ker L_\tau$, so the operator 
$L_\tau$ is Fredholm of index $0$, and therefore self-adjoint by 
Lemma~\ref{lem:selfadj}, for all $\tau\in [0,1]$.
Since the inclusion $E_\tau\into L^2$ is compact, the eigenvectors
$e_n^\tau$ of $L_\tau$ form a Hilbert space basis of $L^2$ for all
$\tau\in [0,1]$.  In particular, for $\tau=0$ we get a unitary
isomorphism
\begin{equation}\label{eq:isom}
  l^2\stackrel{\cong}{\longrightarrow}L^2,\qquad e_n\mapsto e_n^0
\end{equation}
that restricts to an isomorphism $h^1\stackrel{\cong}{\longrightarrow}E_0$.
In the sequel we will use this isomorphism to identify $l^2$ with
$L^2$ and $e_n$ with $e_n^0$, so we can view the $U_\tau$ from
Lemma~\ref{lem:Bernd} as unitary isomorphisms of $L^2$ restricting
to isomorphisms of $E_0$. 

The loop $\{L_\tau\}_{\tau\in [0,1]}$ of self-adjoint Fredholm 
operators is the heart of the counterexample. 
%One is tempted to perform a computation showing that the determinant
%bundle over this loop is nontrivial. 
Unfortunately, these operators have varying domains of
definition. Therefore, we will conjugate 
$\{L_\tau\}_{\tau\in [0,1]}$ with a suitable loop of unitary 
operators to bring all the operators from the desired 
loop on the same footing. Let $\{U_\tau\}_{\tau\in [1,2]}$
be the path of unitary operators in $L^2$ obtained from the one in
Lemma~\ref{lem:Bernd} via identification~\eqref{eq:isom}. 

We define a modified loop $\{\tilde E_\tau\}_{\tau\in [0,2]}$ of subspaces 
of $L^2$ by  
$$
\tilde E_\tau := \begin{cases}
  E_\tau & \tau\in [0,1],\\
  E_0=E_1 & \tau\in [1,2].
\end{cases}
$$
We define a loop $\{\Phi_\tau\}_{\tau\in [0,2]}$ of 
isomorphisms of $L^2$ that restrict to isomorphisms  
$$
  \Phi_\tau|_{E_0}:E_0\longrightarrow \tilde E_\tau
$$
by 
$$
\Phi_\tau f := \begin{cases}
  e^{-i\pi \tau t}f & \tau\in [0,1],\\
  U_\tau f & \tau\in [1,2].
\end{cases}
$$
Note that $e^{-i\pi t}e_n = e^{-i\pi t}e^{i\pi nt} = e_{n-1} = U_1e_n$
and $U_2={\rm Id}$, so $\{\Phi_\tau\}_{\tau\in [0,2]}$ is indeed a
continuous loop. 
Finally, we define 
$$
  \tilde L_\tau:=L|_{\tilde E_\tau}
$$
and 
$$
  T_\tau:=\Phi_\tau^{-1}\tilde L_\tau\Phi_\tau
$$
for $\tau\in [0,2]$. The definition of $\Phi_\tau$
and Lemma~\ref{lem:Bernd} imply that 
$\{T_\tau\}_{\tau\in [0,2]}$ is a loop of self-adjoint operators
in $L^2$ with domain of definition $E_0$. This puts 
us in the setting of Section~\ref{sec:det} with $F:=L^2$,
$E:=E_0$ and $\{T_\tau\}_{\tau\in [0,2]}\subset \FF_s(E,F)$.
Note that 
\begin{equation*}\label{eq:kerTs}
\ker T_\tau=
\begin{cases}
  \R e_0 & \tau=0,\\
  0 & \tau\in (0,1),\\
  \R U_\tau^{-1}e_0 & \tau\in [1,2].
\end{cases}
\end{equation*}
To stabilize the loop $\{T_\tau\}_{\tau\in [0,2]}$, we fix some
$a,b>0$ and consider the vector
$$ 
  G:=ae_0+be_1.
$$
The positivity in Lemma~\ref{lem:Bernd}~(iii)
and the positivity of $a$ and $b$ imply that 
$$
\la G,U_\tau^{-1}E_0\ra>0
$$
for all $\tau\in [1,2]$. Therefore, the vector $G$ is transverse to
$\im T_\tau=(\ker T_\tau)^\perp$ for all $\tau\in [1,2]$, so the
stabilized operator 
$$
  \hat T_\tau:E_0\oplus\R\to L^2,\qquad(f,\zeta)\mapsto 
  T_\tau f+\zeta G
$$
is surjective for all $\tau\in [0,2]$. 
The discussion at the beginning of Section~\ref{sec:det}
implies that the determinant bundle of 
$\{T_\tau\}_{\tau\in [0,2]}$ is isomorphic to the kernel 
bundle of $\{\hat T_\tau\}_{\tau\in [0,2]}$. This loop splits
naturally into two parts. 

{\em The part $\tau\in [0,1]$:}
Using $-i\frac{d}{dt}e_n=\pi n e_n$, we compute for
$f=\sum_{n\in\Z}c_ne_n\in E_0$:
\begin{align*}
  T_\tau f
  &= e^{i\pi n\tau t}\Bigl(-i\ddt\Bigr) e^{-i\pi n\tau t}f
  = e^{i\pi n\tau t}\Bigl(-\pi\tau e^{-i\pi \tau t}f - e^{-i\pi \tau t}i\ddt f\Bigr)\cr
  &= -\pi \tau f -i\ddt f 
  = \sum_{n\in\Z}\pi(n-\tau)c_ne_n. 
\end{align*}
Hence the equation
$$
  \hat T_\tau(f,\zeta) = \sum_{n\in\Z}\pi(n-\tau)c_ne_n +
  \zeta(ae_0+be_1) = 0
$$
is equivalent to $c_n=0$ for all $n\ne 0,1$ and
$$
  \pi(-\tau)c_0+\zeta a = \pi(1-\tau)c_1 + \zeta b = 0. 
$$
Thus $\ker \hat T_\tau$ is described by the equations
$$
  \zeta = \frac{1}{a}\pi\tau c_0 = \frac{1}{b}\pi(\tau-1)c_1
$$
and
$$
  (c_0,c_1)\in \R\bigl(a(\tau-1),b\tau\bigr).
$$
{\em The part $\tau\in [1,2]$:}
Here we write $U_\tau f=\sum_{n\in\Z}c_ne_n$. A computation similar to
the one above shows that 
$$
  \hat T_\tau(f,\zeta) = \sum_{n\in\Z}\pi n c_ne_n +
  \zeta(ae_0+be_1) = 0
$$
is equivalent to $c_n=0$ for all $n\ne 0,1$ and
$$
  0+\zeta a = \pi c_1 + \zeta b = 0. 
$$
The last equations mean that $\zeta=c_1=0$, and therefore
$$
  \ker \hat T_\tau = \R U_\tau^{-1}e_0\oplus 0
  = \R\bigl(c_0(\tau)e_0+c_1(\tau)e_1\bigr)\oplus 0
$$
with coefficients $c_i(\tau)\in \R$ which according to 
Lemma~\ref{lem:Bernd}~(iii) satisfy $c_i(\tau)>0$ for $\tau\in (1,2)$
and $i=0,1$.
%This allows us to summarize the behaviour of $\ker \hat T_\tau$ 
%as follows. The line $\ker \hat T_\tau\subset \R\{e_0,e_1\}$
%interpolates between
%$\R\{e_1\}$ for $\tau=1$ and $\R\{e_0\}$ for $\tau=0$ 
%\begin{center}
%through {\em first and third} quadrants
%\end{center}
%as $\tau$ increases from $1$ to $2$.

After rescaling we may assume $c_1(1)=b$ and $c_0(2)=a$. Then the kernels $\{\ker\hat
T_\tau\}_{\tau\in [0,2]}$ are spanned by the continuous section
$(f(\tau),\zeta(\tau))$ with
$$
f(\tau) = \begin{cases}
  a(\tau-1)e_0 + b\tau e_1 & \tau\in[0,1],\\
  c_0(\tau)e_0+c_1(\tau)e_1 & \tau\in[1,2],
\end{cases} \qquad
\zeta(\tau) = \begin{cases}
  \pi\tau(\tau-1) & \tau\in[0,1],\\
  0 & \tau\in[1,2].
\end{cases}
$$
Since $f(0)=-ae_0$ and $f(2)=ae_0$, this shows nontriviality of the
kernel bundle $\{\ker\hat T_\tau\}_{\tau\in [0,2]}$ over the circle,
and thus of the determinant bundle of $\{T_\tau\}_{\tau\in [0,2]}$. 
In view of isomorphism~\eqref{eq:isom}, this concludes the proof of
Proposition~\ref{prop:counterex}. 
\hfill$\square$

\begin{rem}\label{rem:topology}
For $E,F$ as in Theorem~\ref{thm:det} consider the spaces of operators
\begin{equation}\label{eq:spaces}
   \FF_s^{>\mathfrak{R}}(E,F)^* \subset \FF_s^{>\mathfrak{R}}(E,F)
   \subset \FF_s(E,F) \subset \FF_0(E,F),
\end{equation}
where $\FF_0(E,F)$ denotes the space of Fredholm operators of index zero.
The space $\FF_0(E,F)$ is connected and $\pi_1\FF_0(E,F)=\Z/2\Z$, see
e.g.~\cite{fitzpatrick-pejsachowicz} 
%marginpar{\tiny Check the  reference!}
Proposition~1.3.5 and Corollary~1.5.10.
(here and in the following $\pi_1$ is the fundamental group based at the identity).
Proposition~\ref{prop:counterex} shows that $\pi_1\FF_s(E,F)$ and the
map $\pi_1\FF_s(E,F) \to \pi_1\FF_0(E,F)$ induced by the inclusion are nontrivial.
On the other hand, $\FF_s^{>\mathfrak{R}}(E,F)^*$ has infinitely many
connected components given by the numbers $i(T)$ of negative
eigenvalues, cf.~Remark~\ref{rem:invertible-or}. 
It would be interesting to determine other homotopy groups of the
spaces in~\eqref{eq:spaces}. For example, if one could show
that $\FF_s^{>\mathfrak{R}}(E,F)$ is simply connected, then this would
provide an alternative proof of the orientability of the determinant
bundles in Theorem~\ref{thm:det}. 
\end{rem}

\begin{rem}
The non-orientability phenomenon described in Proposition~\ref{prop:counterex}
arises for example in the algebraic count of intersections of
nonorientable Lagrangians, or of Reeb chords between nonorientable Legendrians. 
\end{rem}

%%%%%%%%%%%%%%%%%%%%%%%%%%%%%%%%%%%%%%%%%%%%%%%%%%%%%%%%%%%%%%%%%%%%%%%%%%
\subsection{The Euler number of gradient vector fields}\label{sec:Euler}
%%%%%%%%%%%%%%%%%%%%%%%%%%%%%%%%%%%%%%%%%%%%%%%%%%%%%%%%%%%%%%%%%%%%%%%%%%

%%%
%\subsection{Gradient vector fields}
%%%

Consider now a Hilbert manifold $X$ (an open subset of a Hilbert space
will suffice for our purposes) and a Hilbert space bundle $E\to X$
with a continuous bundle inclusion $TX\subset E$ such that $T_xX\subset
E_x$ is dense and the inclusion $T_xX\into E_x$ is compact for each
$x\in X$. Let $E$ be equipped with a fibrewise Hilbert
space inner product $\la\ ,\ \ra_E$. In our application, $X$ and $E$
will be maps of Sobolev class $H^2$ and $L^2$, respectively. 

Let $f:X\to\R$ be a $C^1$-map which has an {\em $E$-gradient of class
  $C^1$}, i.e., a $C^1$-section $\Nabla{E}f:X\to E$ satisfying
$$
   Df(x)w = \la\Nabla{E}f(x),w\ra_E \quad \text{for all }w\in T_xX,
   x\in X.
$$
Note that this condition uniquely determines $\Nabla{E}f$. Taking a
derivative of this equation at $x\in\Crit(f)$ we obtain
$$
   D^2f(x)(v,w) = \la D\Nabla{E}f(x)v,w\ra_E \quad \text{for all }v,w\in T_xX,
   x\in\Crit(f),
$$
where $D\Nabla{E}f(x):T_xX\to E_x$ is the linearization of the section
$\Nabla{E}f$ at $x$. Hence the second derivative $D^2f(x)$ is
the composition of the continuous linear maps
$$
\xymatrix
@C=40pt 
{
T_xX\otimes T_xX\into T_xX\otimes E_x
   \ar[r]^{\ \ \ \ \ \ \ \ \ \ \ D\Nabla{E}f(x)\otimes{\rm id}} & E_x\otimes E_x
   \ar[r]^{\ \ \ \ \la\ ,\ \ra_E} & \R. 
}
$$
This implies that $f:X\to\R$ is of class $C^2$, so its second
derivative $D^2f(x)(v,w)$ is symmetric in $v,w$ (this follows from the
usual finite dimensional
result applied to the map $(s,t)\mapsto f\circ\exp_x(sv+tw)$). 
Therefore, the linear map 
$$
D\Nabla{E}f(x):T_xX\to E_x
$$
is symmetric in
the sense of Section~\ref{sec:det}. For $\mathfrak{R}\in\R$ we denote by
$$
   \Func_s^{>\mathfrak{R}}(X) \subset C^2(X,\R)
$$ 
the space of $C^1$-functions $f:X\to\R$ with $E$-gradient of class
$C^1$ such that $D\Nabla{E}f(x)\in\FF_s^{>\mathfrak{R}}(T_xX,E_x)$ for each $x\in\Crit(f)$. 
Then Theorem~\ref{thm:det} implies

\begin{cor}\label{cor:det}
The real line bundle 
$$
\LL\to
\{(f,x)\in\Func_s^{>\mathfrak{R}}(X)\times X\mid Df(x)=0\}
$$ 
associating to $(f,x)$ 
the determinant line $\det(D\Nabla{E}f(x))$ has a canonical orientation. 
\hfill$\square$
\end{cor}

For a critical point $x$ of $f\in\Func_s^{>\mathfrak{R}}(X)$ we denote
by $\ind(x)\in\N_0$ the maximal dimension of a subspace of $T_xX$ on
which $D^2f(x)$ is negative definite, or equivalently, the number of
negative eigenvalues (with multiplicity) of $D\Nabla{E}f(x)$.
%We denote by
%$$
%   \Func_s^+(X)\subset\Func_s(X)
%$$
%the subset of $f$ for which all critical points have finite index. 
A critical point $x$ is called {\em nondegenerate} if
$D\Nabla{E}f(x):T_xX\to E_x$ is invertible. 

\begin{thm}\label{thm:zero-Fredholm}
To each $f\in\Func_s^{>\mathfrak{R}}(X)$ with compact critical point set
we can associate its {\em Euler number} $\chi(\Nabla{E}f)\in\Z$ 
which is uniquely characterized by the following axioms:

(Transversality) 
If all critical points of $f$ are nondegenerate, then there are only
finitely many critical points and
$$
   \chi(\Nabla{E}f) = \sum_{x\in\Crit(f)}(-1)^{\ind(x)}.
$$
(Excision)
For any open neighbourhood $\wt X\subset X$ of $\Crit(f)$ we have
$\chi(\Nabla{E}f) = \chi(\Nabla{E}f|_{\wt X})$. 

%(Cobordism)
%If $W$ is a Hilbert manifold with boundary and $F:W\to Y$ a
%$C^1$-Fredholm map of index $1$ with compact zero set $F^{-1}(0)$,
%then $\chi(F|_{\p W}:\p W\to Y)=0$. 

(Homotopy)
If $F:[0,1]\times X\to\R$ is a $C^1$-map with $E$-gradient of class
$C^1$ and compact critical point set such that
$f_t=F(t,\cdot)\in\Func_s^{>\mathfrak{R}}(X)$ for each $t\in[0,1]$,
then $\chi(\Nabla{E}f_0)=\chi(\Nabla{E}f_1)$. 
\end{thm}

\textbf{Proof: }
The proof consists in applying the proof of~\cite[Theorem
  C.1]{cieliebak-frauenfelder-volkov} to the $E$-gradients, and
upgrading it to integer coefficients using the orientations of the
determinant bundles. Here are the details. 

Consider $f\in\Func_s^{>\mathfrak{R}}(X)$ with compact critical point set. Its
$E$-gradient defines a $C^1$-Fredholm section $\Nabla{E}f:X\to E$ of
index zero with compact zero set $(\Nabla{E}f)^{-1}(0)=\Crit(f)$. 
By Kuiper's theorem, we can trivialize the Hilbert space bundle
$E\cong X\times Y$ and thus view $\Nabla{E}f$ as a $C^1$-Fredholm map
$X\to Y$ to a Hilbert space $Y$. This map satisfies the hypotheses of 
~\cite[Theorem C.1]{cieliebak-frauenfelder-volkov}, so it has a mod
$2$ Euler number uniquely characterized by the analogues of the 
(Transversality), (Excision) and (Homotopy) axioms. By
Corollary~\ref{cor:det}, the determinant line bundle
$\det(D\Nabla{E}f)\to{\rm Crit}(f)$ is canonically oriented. This implies that the
transversely cut out $0$- and $1$-dimensional submanifolds in the proof
of~\cite[Theorem C.1]{cieliebak-frauenfelder-volkov} inherit canonical
orientations and we get a well-defined integer valued Euler number
$\chi(\Nabla{E}f)\in\Z$ (see~\cite{cieliebak-mundet-salamon}). 
The formula in the (Transversality) axiom follows from
Remark~\ref{rem:invertible-or}. 
\hfill $\square$

\begin{rem}
We have stated Corollary~\ref{cor:det} and
Theorem~\ref{thm:zero-Fredholm} only for gradient vector fields
because this covers their applications in this paper. 
Their proofs actually give more general statements for $C^1$-sections
$S:X\to E$ such that $DS(x)\in\FF_s^{>\mathfrak{R}}(T_xX,E_x)$ for
each $x\in S^{-1}(0)$. These should be compared to results on
orientability of Fredholm maps and their mapping degree in the
literature, see
e.g.~\cite{fitzpatrick-pejsachowicz-rabier, abbondandolo-rot} and the
references therein.
\end{rem}

%%%%%%%%%%%%%%%%%%%%%%%%%%%%%%%%%%%%%%%%%%%%%%%%%%%%%%%%%%%%%%%%%%%%%%%%%%
\section{Frozen planet orbits with instantaneous interaction}
%%%%%%%%%%%%%%%%%%%%%%%%%%%%%%%%%%%%%%%%%%%%%%%%%%%%%%%%%%%%%%%%%%%%%%%%%%

%marginpar{NEW}
In this section we consider the real helium atom with the {\em
  instantaneous interaction} between the electrons according to their
Coulomb repulsion. {\em Frozen planet orbits} of this system are smooth maps
$q=(q_1,q_2):S^1\to\R^2$ satisfying
\begin{equation}\label{eq:inst-interaction}
\left\{\;
\begin{aligned}
  \ddot{q}_1(t) &= -\frac{2}{q_1(t)^2}+\frac{1}{(q_1(t)-q_2(t))^2}, \cr
  \ddot{q}_2(t) &= -\frac{2}{q_2(t)^2}-\frac{1}{(q_1(t)-q_2(t))^2}
\end{aligned}
\right.
\end{equation}
as well as the condition
\begin{equation}\label{eq:inst-ineq}
  q_1(t) > q_2(t)\geq 0\qquad\text{for all }t\in S^1. 
\end{equation}
Thus $q_1$ describes the outer electron and $q_2$ the inner electron, 
where the latter undergoes collisions with the nucleus at the origin.
To regularize these collisions, the following setup was introduced
in~\cite{cieliebak-frauenfelder-volkov}.  
One considers the space
\begin{equation}\label{eq:Hin}
\begin{aligned}
  \mathcal{H}_{in} := \Bigl\{
  &z=(z_1,z_2) \in H^2(S^1,\mathbb{R}^2)\;\Bigl|\; 
  \|z_1\|>0,\,\,\|z_2\|>0,\,\, \cr
  &\ \ z_1^2(\tau_{z_1}(t))-z_2^2(\tau_{z_2}(t))>0 \,\,\text{for
    all}\,\,t\in S^1\Bigr\},
\end{aligned}
\end{equation}
where $z_i$ corresponds to the Levi-Civit\`a regularization of $q_i$
and $\tau_{z_i}=t_{z_i}^{-1}$ are the time reparametrizations defined
in Section~\ref{sec:LC}.  
Note that $\mathcal{H}_{in}$ is an open subset of the Hilbert space
$H^2(S^1,\mathbb{R}^2)$ and the last condition in its definition
corresponds to condition~\eqref{eq:inst-ineq}. Integrating this
condition we see that $\mathcal{H}_{in}\subset \mathcal{H}_{av}$.
The {\em instantaneous interaction functional} $\mathcal{B}_{in}
\colon \mathcal{H}_{in} \to \mathbb{R}$ is defined by 
\begin{equation}\label{eq:Bin}
  \mathcal{B}_{in}(z_1,z_2) := \sum_{i=1}^2\Bigl(2\|z_i\|^2\|z_i'\|^2
  + \frac{2}{\|z_i\|^2}\Bigr) 
  - \int_0^1\frac{1}{z_1^2(\tau_{z_1}(t)) - z_2^2(\tau_{z_2}(t))}dt.
\end{equation}
It is proved in~\cite{cieliebak-frauenfelder-volkov} that under the
Levi-Civit\`a transformations $q_i(t)=z_i(\tau_i(t))$ critical points
of $\BB_{in}$ correspond to solutions of~\eqref{eq:inst-interaction}
and~\eqref{eq:inst-ineq}. 

We interpolate between the mean and instantaneous interaction
functionals by 
$$
  \BB_r := r\BB_{in}+
  (1-r)\BB_{av}:\HH_{in}\to\R,\qquad r\in[0,1]. 
$$
It is proved in~\cite{cieliebak-frauenfelder-volkov} that for each
$r\in[0,1]$ the $L^2$-gradient
$$
  \nabla\BB_r = r\nabla\BB_{in}+(1-r)\nabla\BB_{av}:\HH_{in}\to L^2(S^1,\R^2)
$$
is a $C^1$-Fredholm map of index zero. 

As in~\cite[\S6]{cieliebak-frauenfelder-volkov} and
\S\ref{sec:uniqueness} above, we remove the symmetries of $\BB_r$ by
restricting it to a suitable subspace. 
%Following Sections~6.2 and ~6.3 of~\cite{cieliebak-frauenfelder-volkov},
For $k\in\N_0$ we introduce the Hilbert space of {\em symmetric loops} 
\begin{align*}
  H_{\rm sym}^k(S^1,\R^2) 
  := \bigl\{&z=(z_1,z_2)\in H^k(\R/2\Z,\R^2)\;\bigl| z_1(1+\tau) =
  z_1(\tau) = z_1(1-\tau) \\
  &\text{ and }
  -z_2(1+\tau) = z_2(\tau) = z_2(1-\tau) \text{ for all }\tau\bigr\}.
\end{align*}
We consider on $H_{\rm sym}^2(S^1,\R^2)$ the $L^2$-inner product $\la
z,w\ra = \sum_{i=1}^2\int_0^1z_i(\tau)w_i(\tau)d\tau$ and define the
open subset 
\begin{align*}
  X := \bigl\{&z=(z_1,z_2)\in H_{\rm sym}^2(S^1,\R^2) \mid \|z_1\|>0,\,\,\|z_2\|>0,\,\, \cr
  &\ \ z_1^2(\tau_{z_1}(t))-z_2^2(\tau_{z_2}(t))>0 \,\,\text{for all}\,\,t, \cr
  &\ \ z_2'(0)>0,\;z_i(\tau)>0\text{ for all }\tau\in(0,1)\text{ and }i=1,2\bigr\}.
\end{align*}
Note that the first two lines in the definition of $X$ correspond to
the conditions for $\HH_{in}$, and the third line implies that $z_2$
is simple and normalized in the sense of \S\ref{sec:uniqueness}. We
will refer to critical points of $\BB_r$ on $X$ and their Levi-Civit\`a
transforms as {\em normalized simple symmetric frozen planet orbits}
for $\BB_r$.  

It is proved in~\cite{cieliebak-frauenfelder-volkov} that for each
$r\in[0,1]$ the $L^2$-gradient $\nabla\BB_r:X\to H^0_{\rm sym}(S^1,\R^2)$
is a $C^1$-Fredholm map of index zero, and the critical
point set $\mathcal{Z}$ of the $C^2$-map
\begin{equation}\label{eq:BB}
   [0,1]\times X\to\R,\qquad (r,z)\mapsto \BB_r(z)
\end{equation}
is compact. We wish to apply Theorem~\ref{thm:zero-Fredholm} to this
map. For this, we need the following lemma. 

\begin{lemma}\label{lem:spectral-bound}
There exists a constant 
$\mathfrak{R}<0$ such that for each
$(r,z)\in\mathcal{Z}$ the spectrum of the
Hessian $D\nabla\BB_r(z)$ is contained in $(\mathfrak{R},\infty)$.
\end{lemma}

{\bf Proof: }
We claim that the Hessian at
$(r,z)\in\mathcal{Z}$ has the form
$$
  D\nabla\BB_r(z) = \mathfrak{P}_z + \mathfrak{K}_{r,z}: H_{\rm
    sym}^2(S^1,\R^2)\to H_{\rm sym}^0(S^1,\R^2),
$$
where the leading order term is given on $v=(v_1,v_2)\in H_{\rm
    sym}^2(S^1,\R^2)$ by
$$
  \mathfrak{P}_zv = (-4||z_1||^2v_1'',-4||z_2||^2v_2''),
$$
and the lower order term $\mathfrak{K}_{r,z}$ extends to a bounded
operator $H_{\rm sym}^1(S^1,\R^2)\to H_{\rm sym}^0(S^1,\R^2)$
depending continuously on $(r,z)$.
In~\cite{cieliebak-frauenfelder-volkov} it is shown that the lower order term 
is compact, but a more careful inspection of the compactness argument 
reveals that it actually shows the extension to $H^1$. It is enough to see this for $\BB_{av}$ and 
$\BB_{in}$ separately. Below
we abbreviate $H_{\rm sym}^1(S^1,\R^2)$ as $H^1$ etc. 

%\vspace{0.5pt}

For $\BB_{av}$ this follows from 
equations~(93) and~(94) 
of~\cite{cieliebak-frauenfelder-volkov}. 
For $\BB_{in}$ it follows from
the analysis of the Hessian of 
the instantaneous interactions term 
carried out in Sections~A3 to~A.5 of~\cite{cieliebak-frauenfelder-volkov}.
The main additions on top of what is already done
in~\cite{cieliebak-frauenfelder-volkov} are the following:
\begin{itemize}
\item It is written 
in~\cite{cieliebak-frauenfelder-volkov} that the map 
$v\mapsto D_{z_1}v$
defined by formula~(99) is a 
bounded linear operator 
$L^2\mapsto H^1$. This map 
should be used as such when referred
to in Sections~A3 to~A.5
\item The map $v\mapsto Xv$
defined by formula~(109) 
extends as a 
bounded linear operator
$H^1\mapsto H^1$. 
\end{itemize}
This shows the claim above. Now compactness of $\mathcal{Z}$
implies the existence of uniform constants $\delta,C>0$ (independent
of $r,z$) such that
$$
  \la\mathfrak{P}_zv,v\ra
  = \sum_{i=1}^2\la -4\|z_i\|^2v_i'',v_i\ra
  = \sum_{i=1}^24\|z_i\|^2\|v_i'\|^2
  \ge \delta\|v'\|^2
$$
and
$$
  \|\mathfrak{K}_{r,z}v\|\le C \|v\|_{H^1} =  C(\|v\| + \|v'\|).
$$
It follows that
\begin{align*}
  \la D\nabla\BB_r(z)v,v\ra
  &\geq \delta\|v'\|^2 - C(\|v\| + \|v'\|)\|v\| \cr
  &= \delta\Bigl(\|v'\| - \frac{C}{2\delta}\|v\|\Bigr)^2 -
  \Bigl(C+\frac{C^2}{4\delta}\Bigr)\|v\|^2 \cr 
  &\geq - \Bigl(C+\frac{C^2}{4\delta}\Bigr)\|v\|^2,
\end{align*}
and thus
$$
  \la \Bigl(D\nabla\BB_r(z)+C+\frac{C^2}{4\delta}\Bigr)v,v\ra \geq 0.
$$
This shows that the spectrum of $D\nabla\BB_r(z)$ is contained in
$(\mathfrak{R},\infty)$ for any $\mathfrak{R}<-C-\frac{C^2}{4\delta}$. 
\hfill$\square$

By the preceding discussion and Lemma~\ref{lem:spectral-bound}, the
map~\eqref{eq:BB} satisfies the hypotheses of
Theorem~\ref{thm:zero-Fredholm}. It follows that the integral Euler
number $\chi(\nabla\BB_r)\in\Z$ is defined for each $r\in[0,1]$ and
independent of $r$. Now by Corollary~\ref{cor:uniqueness-mean},
the functional 
$\BB_{0}=\BB_{av}$ has a unique
normalized simple symmetric critical point, which is nondegenerate of
index zero.  Therefore, Theorem~\ref{thm:zero-Fredholm} implies the
following result which corresponds to Corollary B from the Introduction.

\begin{cor}\label{cor:frozen-planet-count}
The integral count of normalized simple symmetric frozen planet orbits equals
\begin{equation*}%\label{eq:chi}
   \chi(\nabla\BB_{in})  = \chi(\nabla\BB_{av}) = 1 \in \Z.
\end{equation*}
\hfill$\square$
\end{cor}

\appendix

%%%%%%%%%%%%%%%%%%%%%%%%%%%%%%%%%%%%%%%%%%%%%%%%%%%%%%%%%%%%%%%%
\section{Elliptic integrals}\label{sec:ell-int}
%%%%%%%%%%%%%%%%%%%%%%%%%%%%%%%%%%%%%%%%%%%%%%%%%%%%%%%%%%%%%%%%

For $n \in \mathbb{N}_0$ we consider the elliptic integrals
$$I_n \colon (-\infty,1) \to \mathbb{R}, \quad
m \mapsto
\int_0^1 \frac{\zeta^{2n}}{\sqrt{(1-\zeta^2)(1-m\zeta^2)}}d\zeta.$$
At $m=0$ these integrals can be computed in terms of Euler's beta function and evaluated elementarily,
namely
\begin{eqnarray} \nonumber
I_n(0)&=&\int_0^1 \frac{\zeta^{2n}}{\sqrt{1-\zeta^2}}d\zeta\\ \nonumber
&=&\frac{1}{2}\int_0^1 \xi^{n-\frac{1}{2}}(1-\xi)^{-\frac{1}{2}}d\xi\\ \nonumber
&=&\frac{B\big(n+\frac{1}{2},\frac{1}{2}\big)}{2}\\ \nonumber
&=&\frac{\Gamma(n+\frac{1}{2}\big)\Gamma\big(\frac{1}{2}\big)}{2\Gamma(n+1)}\\ \nonumber
&=&\frac{(2n-1)!! \Gamma\big(\frac{1}{2}\big)^2}{2^{n+1}n!}\\ %\nonumber
&=&\frac{(2n-1)!! \pi}{2^{n+1}n!}.\label{be}
\end{eqnarray}
Here $(2n-1)!!$ equals $(2n-1)(2n-3)\cdots 1$ for $n\geq 1$ and $1$
for $n=0$. 
For $m$ different from zero these elliptic integrals can be expressed
via elliptic integrals of the first and second kind. These are defined
for $m \in (-\infty,1)$ by 
$$K(m)=\int_0^1 \frac{1}{\sqrt{(1-\zeta^2)(1-m\zeta^2)}}d\zeta,
\qquad E(m)=\int_0^1 \frac{\sqrt{1-m\zeta^2}}{\sqrt{1-\zeta^2}}d\zeta.$$
For $n$ equal to zero or one this is obvious. Indeed, we have just
\begin{equation}\label{i0}
I_0(m)=K(m)
\end{equation}
for any $m \in (-\infty,1)$ and if $m \neq 0$, then 
\begin{equation}\label{i1}
I_1(m)=\frac{K(m)-E(m)}{m}.
\end{equation}
For larger $n$ this follows from the recursion formula
\begin{equation}\label{rec}
I_{n+2}(m)=\frac{2(n+1)(m+1)I_{n+1}(m)}{(2n+3)m}-\frac{(2n+1)I_n(m)}{(2n+3)m}
\end{equation}
which allows to express $I_n$ with the help $I_0$ and $I_1$ and consequently in terms
of $K$ and $E$ using (\ref{i0}) and (\ref{i1}). The recursion formula (\ref{rec}) follows from
\begin{eqnarray*}
0&=&\int_0^1 \frac{d}{d\xi}\Big(\xi^{n+\frac{1}{2}}\sqrt{1-(m+1)\xi+m\xi^2}\Big)d\xi\\
&=&\int_0^1 \Bigg( \big(n+\tfrac{1}{2}\big)\xi^{n-\frac{1}{2}}\sqrt{1-(m+1)\xi+m\xi^2}+\frac{\xi^{n+\frac{1}{2}}\big(2m\xi-m-1\big)}{2\sqrt{1-(m+1)\xi+m\xi^2}}
\Bigg)d\xi\\
&=&\int_0^1 \frac{(2n+1)\xi^{n-\frac{1}{2}}\big(1-(m+1)\xi+m\xi^2\big)+\xi^{n+\frac{1}{2}}\big(2m\xi-m-1\big)}{2\sqrt{1-(m+1)\xi+m\xi^2}}d\xi\\
&=&\int_0^1 \frac{(2n+3)m\xi^{n+\frac{3}{2}}-(2n+2)(m+1)\xi^{n+\frac{1}{2}}+(2n+1)\xi^{n-\frac{1}{2}}}{2\sqrt{(1-\xi)(1-m\xi)}}d\xi\\
&=&\int_0^1 \frac{(2n+3)m\zeta^{2n+3}-(2n+2)(m+1)\zeta^{2n+1}+(2n+1)\zeta^{2n-1}}{
\sqrt{(1-\zeta^2)(1-m\zeta^2)}}\zeta d\zeta\\
&=&\int_0^1 \frac{(2n+3)m\zeta^{2(n+2)}-(2n+2)(m+1)\zeta^{2(n+1)}+(2n+1)\zeta^{2n}}{
\sqrt{(1-\zeta^2)(1-m\zeta^2)}}d\zeta\\
&=&(2n+3)mI_{n+2}(m)-2(n+1)(m+1)I_{n+1}(m)+(2n+1)I_n(m).
\end{eqnarray*}
In particular, we obtain
\begin{equation}\label{i2}
I_2(m)=\frac{2(m+1)I_1(m)}{3m}-\frac{I_0(m)}{3m},
\end{equation}
which with the help of (\ref{i0}) and (\ref{i1}) we can alternatively write as
$$I_2(m)=\frac{(m+2)K(m)-2(m+1)E(m)}{3m^2}.$$
The derivatives of the elliptic integrals of the first and second kind can be expressed
as a linear combination of them with $m$-dependent coefficients as follows
\begin{equation}\label{der}
K'(m)=\frac{E(m)-(1-m)K(m)}{2m(1-m)}, \qquad E'(m)=\frac{E(m)-K(m)}{2m}.
\end{equation}
For $E$ this is a straightforward application of (\ref{i1}). Indeed, 
\begin{eqnarray*}
E'(m)&=&-\frac{1}{2}\int_0^1\frac{\zeta^2}{\sqrt{(1-\zeta^2)(1-m\zeta^2)}}d\zeta\\
&=&-\frac{I_1(m)}{2}\\
&=&\frac{E(m)-K(m)}{2m}.
\end{eqnarray*}
For $K$ this is more involved and follows from the following computation
\begin{eqnarray*}
0&=&-m\int_0^1 \frac{d}{d\zeta}\frac{\zeta \sqrt{1-\zeta^2}}{\sqrt{1-m\zeta^2}}d\zeta\\
&=&-m\int_0^1\frac{(1-\zeta^2)(1-m\zeta^2)-\zeta^2(1-m\zeta^2)+m\zeta^2(1-\zeta^2)}
{\sqrt{(1-\zeta^2)(1-m\zeta^2)^3}}d\zeta\\
&=&-m\int_0^1\frac{1-2\zeta^2+m\zeta^4}{\sqrt{(1-\zeta^2)(1-m\zeta^2)^3}}d\zeta\\
&=&
\int_0^1\frac{m(1-m)\zeta^2-(1-m\zeta^2)^2
+(1-m)(1-m\zeta^2)}{\sqrt{(1-\zeta^2)(1-m\zeta^2)^3}}d\zeta\\
&=&2m(1-m)K'(m)-E(m)+(1-m)K(m).
\end{eqnarray*}
It follows from (\ref{der}) that the quotient of the elliptic
integrals $E$ and $K$ satisfies the Riccati differential equation 
\begin{equation}\label{ric} 
\bigg(\frac{E}{K}\bigg)'=-\frac{1}{2m}+\frac{1}{m}\frac{E}{K}-\frac{1}{2m(1-m)}\bigg(\frac{E}{K}\bigg)^2,
\end{equation}
as the following computation shows:
\begin{eqnarray*}
\bigg(\frac{E}{K}\bigg)'&=&\frac{E' K-E K'}{K^2}\\
&=&\frac{1}{2m}\frac{E}{K}-\frac{1}{2m}-\frac{1}{2m(1-m)}\bigg(\frac{E}{K}\bigg)^2+
\frac{1}{2m}\frac{E}{K}\\
&=&-\frac{1}{2m}+\frac{1}{m}\frac{E}{K}-\frac{1}{2m(1-m)}\bigg(\frac{E}{K}\bigg)^2.
\end{eqnarray*}
Combining (\ref{i0}) and (\ref{i1}) we get
$$\frac{I_1}{I_0}=\frac{1}{m}\bigg(1-\frac{E}{K}\bigg).$$
Differentiating this expression and using (\ref{ric}) we compute
\begin{eqnarray*}
\bigg(\frac{I_1}{I_0}\bigg)'&=&\frac{1}{m^2}\frac{E}{K}-\frac{1}{m^2}-\frac{1}{m}\bigg(\frac{E}{K}\bigg)'\\
&=&\frac{1}{m^2}\frac{E}{K}-\frac{1}{m^2}+\frac{1}{2m^2}-\frac{1}{m^2}\frac{E}{K}+
\frac{1}{2m^2(1-m)}\bigg(\frac{E}{K}\bigg)^2\\
&=&-\frac{1}{2m^2}+\frac{1}{2m^2(1-m)}\bigg(1-m\frac{I_1}{I_0}\bigg)^2
\end{eqnarray*}
so that we obtain for the quotient of $I_1$ and $I_0$ the Riccati differential equation
\begin{equation}\label{ric2}
\bigg(\frac{I_1}{I_0}\bigg)'=\frac{1}{2m(1-m)}-\frac{1}{m(1-m)}\frac{I_1}{I_0}+\frac{1}{2(1-m)}
\bigg(\frac{I_1}{I_0}\bigg)^2.
\end{equation}
We end this appendix with a technical lemma about the function
$$F \colon (-\infty,1) \to \mathbb{R}, \quad m \mapsto (2-m)\frac{I_1}{I_0}(m)$$
which we need to prove our nondegeneracy of critical points of collision type for the frozen
functional for positive parameters $r$. From (\ref{be}) we have
$$I_0(0)=\frac{\pi}{2}, \qquad I_1(0)=\frac{\pi}{4},$$
and therefore
$$F(0)=1.$$
Our technical lemma is the following.
\begin{lemma}\label{mono}
For $m<0$ we have $F(m)>1$. 
\end{lemma}
\textbf{Proof: } We consider the function
$$G \colon (-\infty,1) \to \mathbb{R}, \quad m \mapsto (2-m)I_1(m)$$
and show that it is strictly decreasing for negative $m$. For that purpose we 
compute its derivative
\begin{eqnarray*}
G'(m)&=&-\int_0^1 \frac{\zeta^{2}}{\sqrt{(1-\zeta^2)(1-m\zeta^2)}}d\zeta\\
& &+\frac{2-m}{2}\int_0^1 \frac{\zeta^{4}}{\sqrt{(1-\zeta^2)(1-m\zeta^2)^3}}d\zeta\\
&=&\frac{1}{2}\int_0^1 \frac{(2-m)\zeta^{4}-2\zeta^2(1-m\zeta^2)}{\sqrt{(1-\zeta^2)(1-m\zeta^2)^3}}d\zeta\\
&=&\frac{1}{2}\int_0^1 \frac{\zeta^2\big((2+m)\zeta^2-2\big)}{\sqrt{(1-\zeta^2)(1-m\zeta^2)^3}}d\zeta.
\end{eqnarray*}
If $m$ is negative the enumerator is nonpositive, and strictly negative for
$\zeta \in (0,1)$. This shows that
$$G'(m)<0\ \text{for}\ m<0.$$
Since $I_0=K$ is strictly increasing, we see that
$$F=\frac{G}{I_0}$$
as the quotient of a positive strictly decreasing function and a
positive strictly increasing function is strictly decreasing for
negative values of $m$. Since $F(0)=1$, we conclude that
$$F(m)>1\ \text{for}\ m<0,$$
which proves the lemma. \hfill $\square$

Data availability statement is not applicable.

The authors declare that they do not have any conflict 
of interest.

\end{document}